\documentclass[11 pt,reqno]{amsart}
\usepackage{amssymb}
\usepackage{amsmath}
\usepackage{amsthm}
\usepackage{amscd}
\usepackage{xypic}
\xyoption{all}
\usepackage[usenames,dvipsnames]{color}
\usepackage{hyperref}

\numberwithin{equation}{section}
\newtheorem*{theorem}{Theorem}
\newtheorem*{corollary}{Corollary}
\newtheorem*{lemma}{Lemma}
\newtheorem*{proposition}{Proposition}
\newtheorem*{conjecture}{Conjecture}

\theoremstyle{definition}
\newtheorem*{definition}{Definition}
\newtheorem*{remark}{Remark}

\renewcommand{\proof}{\noindent{\sc Proof. }}
\newcommand {\halmos}{\hfill{$\blacksquare$}}

\newcommand {\g}{\mathfrak g}
\newcommand {\h}{\mathfrak h}
\renewcommand {\u}{\mathfrak u}

\newcommand {\IA}{\mathbb A}
\newcommand {\IC}{\mathbb C}
\newcommand {\IN}{\mathbb N}
\newcommand {\IV}{\mathbb V}
\newcommand {\IR}{\mathbb R}

\newcommand {\IZ}{\mathbb Z}

\newcommand {\A}{\mathcal A}
\newcommand {\D}{\mathcal D}

\renewcommand {\H}{\mathcal H}
\newcommand {\R}{\mathcal R}
\newcommand {\T}{\mathcal T}
\newcommand {\V}{\mathcal V}
\newcommand {\Y}{\mathcal Y}

\newcommand {\sfA}{{\mathsf A}}
\newcommand {\sfB}{{\mathsf B}}

\newcommand {\sfG}{{\mathsf G}}

\newcommand {\bfI}{\mathbf I}
\newcommand {\bfA}{\mathbf A}
\newcommand {\bfD}{\mathbf D}

\newcommand {\ol}[1]{\overline{#1}}
\newcommand {\ul}[1]{\underline{#1}}
\newcommand {\wh}[1]{\widehat{#1}}
\newcommand {\wt}[1]{\widetilde{#1}}

\newcommand {\cbc}{case--by--case }
\newcommand {\fd}{finite--dimensional }
\newcommand {\eg}{{\it e.g. }}
\newcommand {\etal}{{\it et al. }}

\newcommand {\lhs}{left--hand side }
\newcommand {\rhs}{right--hand side }
\newcommand {\wrt}{with respect to }

\newcommand {\FR}{Frenkel--Reshetikhin }

\newcommand {\KZ}{KZ }
\newcommand {\qW}{quantum Weyl group }
\newcommand {\TV}{Tarasov--Varchenko }
\newcommand {\QYBE}{quantum Yang--Baxter equations }
\newcommand {\daha}{degenerate affine Hecke algebra }

\newcommand {\ad}{\operatorname{ad}}
\newcommand {\Ad}{\operatorname{Ad}}
\newcommand {\Aut}{\operatorname{Aut}}

\newcommand {\End}{\operatorname{End}}
\newcommand {\Hom}{\operatorname{Hom}}
\newcommand {\Ker}{\operatorname{Ker}}

\newcommand{\ev}{\operatorname{ev}}
\newcommand {\id}{\operatorname{id}}
\newcommand{\tr}{\operatorname{tr}}
\newcommand{\ind}{\operatorname{ind}}
\newcommand {\rk}{\operatorname{rank}}

\renewcommand {\c}{\mathfrak c}
\newcommand {\Wa}{W_a}
\newcommand {\Ba}{B_a}
\newcommand {\wts}{\wt{s}}
\newcommand {\Qv}{Q^\vee}

\newcommand {\red}[1]{{\widetilde{#1}}^{\operatorname{\scriptscriptstyle{red}}}}

\newcommand {\nablac}{\nabla_{\scriptstyle{\kappa}}}

\newcommand {\reg}{_{\scriptstyle{\operatorname{reg}}}}
\newcommand {\hreg}{\h\reg}

\newcommand {\gl}[1]{\mathfrak{gl}_{#1}}
\renewcommand {\sl}[1]{\mathfrak{sl}_{#1}}

\newcommand {\Uhghat}{U_\hbar(L\g)}

\newcommand {\fml}{[\![\hbar]\!]}
\newcommand {\fmlu}{[\![u^{-1}]\!]}
\newcommand {\ICh}{\IC\fml}
\newcommand {\Ch}{\IC}

\newcommand {\Omit}[1]{}

\newcommand {\res}{\operatorname{res}}

\newcommand {\half}[1]{\frac{#1}{2}}

\newcommand {\aand}{\qquad\text{and}\qquad}

\renewcommand {\SS}{\mathfrak S}

\newcommand {\<}{\langle}
\renewcommand {\>}{\rangle}

\renewcommand {\root}[1]{\alpha_#1}
\newcommand {\cor}[1]{\alpha_{#1}^\vee}
\newcommand {\cow}[1]{\lambda_{#1}^\vee}

\newcommand {\etaa}[1]{\frac{d#1}{e^{#1}-1}}

\newcommand {\etaeta}[2]{\eta_{#1}\wedge\eta_{#2}\,[t_{#1},t_{#2}]}

\newcommand {\height}{\operatorname{ht}}
\newcommand {\sign}{\operatorname{sign}}

\newcommand {\ih}{^\hbar_i}

\newcommand {\TT}{\not{T}}

\begin{document}

\title[The Trigonometric connection of a simple Lie algebra]
{The Trigonometric Casimir connection of a simple Lie algebra}
\author[V. Toledano Laredo]{Valerio Toledano Laredo}
\address
{Department of Mathematics,
Northeastern University,
567 Lake Hall,
360 Huntington Avenue,
Boston,
MA 02115.}
\email{V.ToledanoLaredo@neu.edu}
\dedicatory{To Corrado De Concini, on his 60th birthday}
\thanks{Supported in part by NSF grants DMS--0707212 and DMS--0854792}
\begin{abstract}
Let $\g$ be a complex, semisimple Lie algebra, $G$ the corresponding
simply--connected Lie group and $H\subset G$ a maximal torus. We
construct a flat connection on $H$ with logarithmic singularities on the
root hypertori and values in the Yangian $Y(\g)$ of $\g$. By analogy with
the rational Casimir connection of $\g$, we conjecture that the monodromy
of this trigonometric connection is described by the quantum Weyl group
operators of the quantum loop algebra $U_\hbar(L\g)$.
\end{abstract}
\maketitle

\setcounter{tocdepth}{1}
\tableofcontents

\section{Introduction}

\subsection{} 

Let $\g$ be a complex, simple Lie algebra, $\h\subset\g$ a
Cartan subalgebra, $\Phi\subset\h^*$ the corresponding
root system and $W$ its Weyl group. For
each $\alpha\in\Phi$, let $\sl{2}^\alpha=\<e_\alpha,f_\alpha,
h_\alpha\>\subset\g$ be the corresponding three--dimensional
subalgebra and denote by
$$\kappa_\alpha=\frac{(\alpha,\alpha)}{2}
\left(e_\alpha f_\alpha+f_\alpha e_\alpha\right)$$
its truncated Casimir operator \wrt the restriction to $\sl{2}^
\alpha$ of a fixed non--degenerate, ad--invariant bilinear form
$(\cdot,\cdot)$ on $\g$. Let
$$\hreg=\h\setminus\bigcup_{\alpha\in\Phi}\Ker(\alpha)$$
be the set of regular elements in $\h$, $V$ a \fd $\g$--module,
and $\IV$ the holomorphically trivial vector bundle over $\hreg$
with fibre $V$. Recall that the Casimir connection of $\g$ is
the holomorphic connection on $\IV$ given by
\begin{equation}\label{eq:rational}
\nablac=
d-\hbar\sum_{\alpha\in\Phi_+}\frac{d\alpha}{\alpha}\,\kappa_\alpha
\end{equation}
where $d$ is the de Rham differential, $\Phi_+\subset\Phi$ is a system
of positive roots and $\hbar$ is a complex number. This connection
was discovered independently by C. De Concini around 1995 (unpublished),
by J. Millson and the author \cite{MTL,TL1} and Felder \etal \cite{FMTV},
and shown to be flat for any $\hbar\in\IC$.

The monodromy of $\nablac$ gives representations of the generalised
braid group  $B=\pi_1(\hreg/W)$ which are described by the quantum
Weyl group operators of the quantum group $U_\hbar\g$, a fact which
was conjectured by De Concini (unpublished) and independently in \cite
{TL1,TL2} and proved in \cite{TL1,TL3}.

\subsection{}	

Let $P\subset\h^*$ be the weight lattice of $\g$ and $H=\Hom_\IZ(P,
\IC^*)$ the dual algebraic torus with Lie algebra $\h$ and coordinate
ring given by the group algebra $\IC P$. We denote the function
corresponding to $\lambda\in P$ by $e^\lambda\in\IC[H]$. The main
goal of the present paper is to define a {\it trigonometric} version of
the connection \eqref{eq:rational}, that is a connection defined on
$H$ with the logarithmic forms $d\alpha/\alpha$ replaced by $d\alpha
/(e^\alpha-1)$.

As is well known from the study of Cherednik's affine KZ (AKZ)
connection (see, \eg\cite{Ch3}), both the flatness and $W$--equivariance
of such a connection require that it possess a 'tail', that is be of the
form
\begin{equation}\label{eq:trigo form}
\wh{\nabla}_\kappa=
d-
\hbar\sum_{\alpha\in\Phi_+}\frac{d\alpha}{e^\alpha-1}\,\kappa_\alpha-A
\end{equation}
where $A$ is a translation--invariant one--form on $H$. The analogy
with the AKZ equations further suggests that $A$ should take values
in a suitable extension of the enveloping algebra $U\g$, which is to
$U\g$ what the \daha $\H'$ is to the group algebra $\IC W$.

\subsection{} 

The correct extension turns out to be the {\it Yangian} $Y(\g)$, which
is a deformation of the enveloping algebra $U(\g[t])$ over the ring
$\IC[\hbar]$. Let $\nu:\h\to\h^*$ be the isomorphism determined by the
inner product $(\cdot,\cdot)$, set $t_i=\nu^{-1}(\alpha_i)$, where $\alpha_1,
\ldots,\alpha_n$ are the simple roots of $\g$ relative to $\Phi_+$ and
let $t^i=\cow{i}$ be the dual basis of $\h$ given by the fundamental
coweights. Let $T(u)_r$, $u\in\h$, $r\in\IN$ be the Cartan loop generators
of $Y(\g)$ in Drinfeld's new realisation (see \cite{Dr2} and \S   \ref{se:trigo g}
for definitions). Let $\{u^i\}$ be a basis of $\h$, $\{u_i\}$ the dual basis
of $\h^*$ and regard the differentials $du_i$ as translation--invariant
one--forms on $H$. The main result of this paper is the following

\begin{theorem}
The $Y(\g)$--valued connection on $H$ given by
$$\wh{\nabla}_\kappa=
d-
\hbar\sum_{\alpha\in\Phi_+}\frac{d\alpha}{e^\alpha-1}\,\kappa_\alpha
+du_i\,\left(2T(u^i)_1-\half{\hbar}(u^i,t^j)t_j^2\right)
$$
is flat and $W$--equivariant\footnote{we follow the
standard convention that in any expression involving $u_i$ and $u^i$, or
$t^i$ and $t_i$, summation over $i$ is implicit.}.
\end{theorem}

\subsection{}

We call $\wh{\nabla}_\kappa$ the {\it trigonometric Casimir connection}
of $\g$. Its monodromy defines representations of the {\it affine braid group}
$\wh{B}=\pi_1(H\reg/W)$ on any \fd module over $Y(\g)$, where
$$H\reg=H\setminus\bigcup_{\alpha\in\Phi}\{e^\alpha=1\}$$
is the set of regular elements in $H$. 
By analogy with the rational case, we conjecture that these representations
are equivalent to the quantum Weyl group action of $\wh{B}$ on \fd modules
over the quantum loop algebra $U_\hbar(\g[z,z^{-1}])$.

\subsection{}	

We turn now to a detailed description of the paper.

In Section \ref{se:trigo root}  we obtain a necessary and sufficient condition
for a connection of the form \eqref{eq:trigo form} to be flat and equivariant
under the Weyl group $W$  (Theorem \ref{th:trigo flat} and Proposition
\ref{pr:equivariance}), thus effectively determining the Lie algebras of the
fundamental groups $\pi_1(H\reg)$ and $\pi_1(H\reg/W)$. Consistently
with Cherednik's study of the AKZ connection, the quadratic equations
giving the flatness and equivariance of $\wh{\nabla}_\kappa$ specialise,
when $\kappa_\alpha$ is replaced by the orthogonal reflection $s_\alpha
\in W$ to the defining relations of the \daha corresponding to $W$.

In Section \ref{se:trigo g} we review Drinfeld's two presentations of the
Yangian $Y(\g)$ of a simple Lie algebra $\g$. We then use their interplay
to solve the above quadratic equations in $Y(\g)$, thereby obtaining the
trigonometric Casimir connection of $\g$ (Theorem \ref{th:trigo Casimir}).

In Section \ref{se:monodromy}, we explain how to define monodromy
representations of the affine braid group from $\wh{\nabla}_\kappa$
and conjecture that these are described by the quantum Weyl group
operators of the quantum loop algebra $U_\hbar(\g[z,z^{-1}])$.

In Section \ref{se:trigo gl}, we define a trigonometric connection with
values in the Yangian of $\gl{n}$ by using the interplay between its
loop and RTT presentations (Theorem \ref{th:trigo gl}). We then relate
it to the trigonometric Casimir connection of $\sl{n}$. We also check
that, when computed in a tensor product of $m$ evaluation modules,
it coincides with the trigonometric dynamical differential equations
\cite{TV} which are differential equations on $(\IC^\times)^n$ with
values in $U\gl{n}^{\otimes m}$.

In Section \ref{se:bispesctral}, we show that the trigonometric Casimir
connection commutes with $q$KZ difference equations of \FR determined
by the rational $R$--matrix of $Y(\g)$ (Theorem \ref{th:bispectral}),
a fact which was checked in \cite{TV} for the trigonometric dynamical
differential equations.

In Section \ref{se:AKZ}, we review the definition of the \daha
$\H'$ of $W$ \cite{Lu} and of the corresponding $\H'$--valued
AKZ connection \cite{Ch3}. We then show that if $V$ is a $Y(\g)
$--module whose restriction to $\g$ is {\it small}, that is such that
$2\alpha$ is not a weight for any root $\alpha$ \cite{Br,Re}, the
zero weight space $V[0]$ carries a natural action of $\H'$. Moreover,
the trigonometric Casimir connection with coefficients in $V[0]$
coincides with the AKZ connection with values in this $\H'$--module
(Theorem \ref{th:trigo=AKZ}).

The final appendix, Section \ref{se:Tits} contains a discussion
of the Tits extensions of affine Weyl groups which is needed
for Section \ref{se:monodromy}.

\section{The trigonometric connection of a root system}\label{se:trigo root}

\subsection{General form}

Let $E$ be a Euclidean vector space, $\Phi\subset E^*$ a reduced,
crystallographic root system. Let $Q^\vee\subset E$ be the lattice
generated by the coroots $\alpha^\vee$, $\alpha\in\Phi$ and $P\subset
E^*$ the dual weight lattice. Let $H=\Hom_\IZ(P,\IC^*)$ be the
complex algebraic torus with Lie algebra $\h=\Hom_\IZ(P,\IC)$
and coordinate ring given by the group algebra $\IC P$. We denote
the function corresponding to $\lambda\in P$ by $e^\lambda\in\IC[H]$
and set
\begin{equation}\label{eq:Hreg}
H\reg=H\setminus\bigcup_{\alpha\in\Phi}\{e^\alpha=1\}
\end{equation}

Let $A$ be an algebra endowed with the following data:
\begin{itemize}
\item a set of elements $\{t_\alpha\}_{\alpha\in\Phi}\subset A$ such
that $t_{-\alpha}=t_{\alpha}$
\item a linear map $\tau:\h\to A$
\end{itemize}
Consider the $A$--valued connection on $H\reg$ given by
\begin{equation}\label{eq:trigo}
\nabla=d-
\sum_{\alpha\in\Phi_+}\etaa{\alpha}\,t_{\alpha}-
du_i\,\tau(u^i)
\end{equation}
where $\Phi_+\subset\Phi$ is a chosen system of positive roots,
$\{u_i\}$ and $\{u^i\}$ are dual bases of $\h^*$ and $\h$ respectively,
the differentials $du_i$ are regarded as translation--invariant
one--forms on $H$ and the summation over $i$ is implicit.

\subsection{Positive roots}\label{ss:chamber}

The form of the connection \eqref{eq:trigo} depends upon the
choice of the system of positive roots $\Phi_+\subset\Phi$. Let
however $\Phi'_+\subset\Phi$ be another such system, then

\begin{proposition}\label{pr:chamber}
The connection \eqref{eq:trigo} may be rewritten as
$$\nabla=
d-\sum_{\alpha\in\Phi_+'}\etaa{\alpha}\,t_\alpha-du_i\,\tau'(u^i)$$
where $\tau':\h\to A$ is given by
\begin{equation}
\tau'(v)=\tau(v)-\sum_{\alpha\in\Phi_+\cap \Phi_-'}\alpha(v)\,t_\alpha
\end{equation}
\end{proposition}
\proof Write the second summand in \eqref{eq:trigo} as
$$\sum_{\alpha\in\Phi_+\cap \Phi_+'}\etaa{\alpha}\,t_\alpha-
\sum_{\alpha\in\Phi_-\cap \Phi_+'}\frac{d\alpha}{e^{-\alpha}-1}\,t_{-\alpha}$$
where $\Phi_-=-\Phi_+$. Since
\begin{equation}\label{eq:trick}
\frac{1}{1-e^{-\alpha}}=
\frac{e^\alpha}{e^\alpha-1}=
\frac{1}{e^\alpha-1}+1
\end{equation}
and $t_{-\alpha}=t_\alpha$, the above is equal to
$$\sum_{\alpha\in\Phi_+'}\etaa{\alpha}\,t_\alpha+
\sum_{\alpha\in\Phi_-\cap \Phi_+'}d\alpha\,t_\alpha$$
which yields the required result since $\alpha=u_i\alpha(u^i)$. \halmos

\begin{definition}
If $W$ is the Weyl group of $\Phi$ and $w\in W$ the unique element
such that $\Phi_+'=w\Phi_+$, we denote $\tau'$ by $\tau_w$. Thus,
\begin{equation}\label{eq:chamber}
\tau_w(v)=\tau(v)-\sum_{\alpha\in\Phi_+\cap w\Phi_-}\alpha(v)\,t_\alpha
\end{equation}
\end{definition}

\subsection{Delta form}\label{ss:delta} 

Choose $\Phi_+'=\Phi_-$ in Proposition \ref{pr:chamber}. Comparing
the corresponding expressions for $\nabla$ shows that it may be more
invariantly rewritten as
$$\nabla=
d-\frac{1}{2}\sum_{\alpha\in\Phi}\etaa{\alpha}\,t_\alpha-du_i\,\delta(u^i)$$
where $\delta:\h\to A$ is given by
\begin{equation}\label{eq:delta}
\delta(v)=\tau(v)-\half{1}\sum_{\alpha\in\Phi_+}\alpha(v)\,t_\alpha
\end{equation}
Alternatively, substituting \eqref{eq:delta} into \eqref{eq:trigo} yields
\begin{equation}\label{eq:delta form}
\nabla
=
d-\half{1}\sum_{\alpha\in\Phi_+}
\frac{e^\alpha+1}{e^\alpha-1}d\alpha\,t_\alpha-du_i\,\delta(u^i)
\end{equation}
We shall occasionally refer to \eqref{eq:trigo} and \eqref{eq:delta form}
as the {\it $\tau$ and $\delta$--forms} of the connection $\nabla$
respectively. Note that the latter does not depend upon the choice
of $\Phi_+\subset\Phi$.

\subsection{Root subsystems}

For a subset $\Psi\subset\Phi$ and subring $R\subset\IR$, let $\<\Psi\>
_R\subset E^*$ be the $R$--span of $\Psi$.

\begin{definition}
A {\it root subsystem} of $\Phi$ is a subset $\Psi\subset\Phi$ such that
$\<\Psi\>_\IZ\cap\Phi=\Psi$. $\Psi$ is {\it complete} if $\<\Psi\>_{\IR}\cap
\Phi=\Psi$. If $\Psi\subset\Phi$ is a root subsystem, we set $\Psi_+=\Psi
\cap\Phi_+$.
\end{definition}

\begin{remark} According to the above definition, the short roots of the
root system $\sfB_2$ (resp. $\sfG_2$) are not a root subsystem, but the
long ones constitute a root subsystem of type $\sfA_1\times\sfA_1$ (resp.
$\sfA_2$) which is not complete. Another root subsystem of $\Phi=\sfG_2$
is given by $\{\pm\alpha,\pm\beta\}$ where $\alpha,\beta$ are two orthogonal
roots (necessarily of different lengths).
\end{remark}

\subsection{Integrability}

The following is the main result of this section.

\begin{theorem}\label{th:trigo flat}\hfill
\begin{enumerate}
\item The connection $\nabla$ is flat if, and only if the following relations hold\\

\begin{itemize}
\item For any rank 2 root subsystem $\Psi\subset\Phi$ and $\alpha\in\Psi$,
\begin{equation}
[t_\alpha,\sum_{\beta\in\Psi_+}t_\beta]=0
\tag{$tt$}
\label{eq:quadratic}
\end{equation}
\item For any $u,v\in \h$,
\begin{equation}
[\tau(u),\tau(v)]=0
\tag{$\tau\tau$}
\label{eq:comm}
\end{equation}
\item For any $\alpha\in\Phi_+$, $w\in W$ such that $w^{-1}\alpha$ is a
simple root and $u\in \h$ such that $\alpha(u)=0$,
\begin{equation}
[t_\alpha,\tau_w(u)]=0
\tag{$t\tau$}
\label{eq:cross}
\end{equation}
\end{itemize}
\item Modulo the relations \eqref{eq:quadratic}, the relations \eqref
{eq:cross} are equivalent to
\begin{equation}\label{eq:easy cross}
[t_\alpha,\delta(v)]=0
\tag{$t\delta$}
\end{equation}
for any $\alpha\in\Phi$ and $v\in \h$ such that $\alpha(v)=0$, where
$\delta:\h\to A$ is given by \eqref{eq:delta}.
\end{enumerate}
\end{theorem}
\noindent
The proof of Theorem \ref{th:trigo flat} occupies the paragraphs
\ref{ss:start flat}--\ref{ss:end flat}.

\subsection{} 

We spell out below the relations \eqref{eq:quadratic} in the case
when $\Phi$ is of rank 2. For $\Psi=\Phi$, they read
\begin{equation}\label{eq:Phi}
[t_\alpha,\sum_{\beta\in\Phi_+}t_\beta]=0
\quad\text{for any $\alpha\in\Phi$}
\tag{$\mathsf{\Phi}$}
\end{equation}
In particular, if $\Phi=\{\pm\alpha,\pm\beta\}$ is of type $\sfA_1\times
\sfA_1$ then
\begin{equation}\label{eq:A1A1}
[t_\alpha,t_\beta]=0 \tag{$\sfA_1\times\sfA_1$}
\end{equation}
For $\Phi=\sfB_2$, the long roots $\{\pm\beta_1,\pm\beta_2\}$
form an $\sfA_1\times\sfA_1$ subsystem so that
\begin{equation}\label{eq:A1A1B2}
[t_{\beta_1},t_{\beta_2}]=0
\tag{$\sfA_1\times\sfA_1\subset\sfB_2$}
\end{equation}
For $\Phi=\sfG_2$, there are two types of root subsystems: that
formed by the $\sfA_2$ configuration of long roots $\{\pm\beta_1,
\pm\beta_2,\pm\beta_3\}$, leading to
\begin{equation}\label{eq:A2G2}
[t_{\beta_i},t_{\beta_j}+t_{\beta_k}]=0
\quad\text{for any $i,j,k\in\{1,2,3\}$ distinct}
\tag{$\sfA_2\subset G_2$}
\end{equation}
and the $\sfA_1$ configurations $\{\pm\beta,\pm\gamma\}$ formed
by a long root and an orthogonal short one, leading to
\begin{equation}\label{eq:A1A1G2}
[t_\beta,t_\gamma]=0
\tag{$\sfA_1\times\sfA_1\subset\sfG_2$}
\end{equation}
Combining relations \eqref{eq:Phi}, \eqref{eq:A2G2} and \eqref{eq:A1A1G2}
yields in particular the following relations
\begin{equation}\label{eq:variant}
[t_\beta,t_{\gamma'}+t_{\gamma''}]=0
\end{equation}
where $\beta\in\sfG_2$ is long and $\gamma',\gamma''$ are the short
positive roots which are not orthogonal to $\beta$.

\begin{remark} If $\Phi$ is not simply--laced, the relations \eqref{eq:quadratic}
are stronger than those yielding the flatness of the rational connection
$$\nabla=d-\sum_{\alpha\in\Phi_+}\frac{d\alpha}{\alpha}\,t_\alpha$$
Indeed, the latter involve two dimensional subspaces of $\h^*$
spanned by elements of $\Phi$ \cite{Ko} and therefore only
those rank 2 subsystems of $\Phi$ which are complete. The
relevance of additional relations corresponding to non--complete
subsystems was first pointed out in the closely related context
of the Yang--Baxter equations by Cherednik \cite[\S   6.1]{Ch3}.
\end{remark}

\subsection{} \label{ss:start flat}

Let $Q\subset P$ be the root lattice generated by $\Phi$ and $T
=\Hom_\IZ(Q,\IC^*)$ the corresponding complex algebraic torus
of adjoint type. $T$ has coordinate ring $\IC Q$ and is biregularly
isomorphic to the standard torus $(\IC^*)^n$ by sending $p\in T$
to the point with coordinates $z_i=e^{-\alpha_i}(p)$, where $\alpha_i$
varies over the simple roots of $\Phi$ relative to $\Phi_+$.

\subsection{}\label{ss:connected}

The torus $T$ is a quotient of $H$ and the form of $\nabla$ shows
that it may be regarded as a connection on the trivial vector bundle
with fibre $A$ over $T$. As such, $\nabla$ has singularities on the
codimension one subtori
$$T_\alpha=\{e^\alpha=1\}\subset T$$
where $\alpha\in\Phi$. Given a subset
$\Psi\subset\Phi$, we shall be interested in the connectedness of the
intersection $\bigcap_{\alpha\in\Psi} T_\alpha$. Let $\<\Psi\>_{\IZ}\subset
Q$ be the $\IZ$--span of $\Psi$ and set \cite[\S  3.1]{DCP}
$$\ol{\<\Psi\>}_{\IZ}=
\{\gamma\in Q|\,m\gamma\in\<\Psi\>_{\IZ}\,\,\text{for some $m\in\IZ^*$}\}$$
Since
$$\IC[\bigcap_{\alpha\in\Psi}T_\alpha]=
\IC\,Q/\<\Psi\>_{\IZ}\cong
\IC\,Q/\ol{\<\Psi\>}_{\IZ}\otimes \IC\,\ol{\<\Psi\>}_{\IZ}/\<\Psi\>_{\IZ}$$
the connected components of $\bigcap_{\alpha\in\Psi}T_\alpha$ are tori
labelled by the characters of the finite abelian group $\ol{\<\Psi\>}_{\IZ}/
\<\Psi\>_{\IZ}$. In particular, if $\Psi=\{\alpha\}$, we see that each $T_
\alpha$ is connected since $\alpha$ is indivisible in $Q$.

\subsection{} 

The necessity of relations \eqref{eq:quadratic}--\eqref{eq:cross} follows
from the computation of the residues of the curvature $\Omega$ of
$\nabla$ to be carried out in \S \ref{ss:start nec}--\S \ref{ss:end nec}.

Specifically, write $\nabla=d-A$. Since $dA=0$, $\Omega$ is equal to
$A\wedge A=\Omega_1+\Omega_2+\Omega_3$, where
\begin{align}
\Omega_1&=
\half{1}
\sum_{\alpha,\beta}
\frac{d\alpha}{e^\alpha-1}\wedge\frac{d\beta}{e^\beta-1}\,[t_\alpha,t_\beta]\\
\Omega_2&=
\sum_{\alpha,i}\frac{d\alpha}{e^\alpha-1}\wedge du_i\,[t_\alpha,\tau(u^i)]\\
\Omega_3&=
\half{1}
\sum_{i,j}du_i\wedge du_j\,[\tau(u^i),\tau(u^j)]
\end{align}

\subsection{}\label{ss:start nec}

Let $\alpha\in\Phi$ and denote the inclusion $T_\alpha\hookrightarrow T$
by $\imath_\alpha$. Then
\begin{align*}
\res_{T_\alpha}\Omega_1
&=\imath_\alpha^*\sum_{\beta\neq\alpha}\etaa{\beta}\,[t_\alpha,t_\beta]
\\
\res_{T_\alpha}\Omega_2
&=\imath_\alpha^*du_i\,[t_\alpha,\tau(u^i)]
\end{align*}
and $\res_{T_\alpha}\Omega_3=0$ since $\Omega_3$ is regular on
$T_\alpha$.

\subsection{}

Let $\Psi\subset\Phi$ be a rank 2 root subsystem and set
$$T_\Psi=\bigcap_{\beta\in\Psi}T_\beta$$
By \S \ref{ss:connected}, $T_\Psi$ is a codimension two subtorus of $T$
with group of components $\Hom(\ol{\<\Psi\>}_{\IZ}/\<\Psi\>_{\IZ},\IC^*)$.

Since $\ol{\<\Psi\>}_{\IZ}/\<\Psi\>_{\IZ}$ is cyclic (of order 1, 2 or 3
depending on the type of $\Psi$ and $\ol{\<\Psi\>}_{\IZ}\cap\Phi$),
there exists a character $\chi$ of $\ol{\<\Psi\>}_{\IZ}/\<\Psi\>_{\IZ}$
with trivial kernel. It follows that the corresponding component $T_
\Psi^\chi$ of $T_\Psi$ is contained in some $T_\gamma$ if, and only
if, $\gamma\in\Psi$.

Together with \S \ref{ss:start nec}, this implies that for any $\alpha
\in\Psi$,
$$\res_{T_\Psi^\gamma}\res_{T_\alpha}\Omega=
[t_\alpha,\sum_{\beta\in\Psi_+}t_\beta]$$
thus showing the necessity of \eqref{eq:quadratic}.

\subsection{}\label{ss:compactification}

Let $\ol{T}\cong\IC^n$ be the partial compactification determined by the
embedding $T\hookrightarrow(\IC^*)^n$ given by sending $p\in T$ to the
point with coordinates $z_i=e^{-\alpha_i}(p)$. We wish to determine the
residues of $\Omega$ on the divisors
$$T_i=\{z_i=0\}\subset\ol{T}$$
To this end, we first rewrite $\nabla$ in the coordinates $z_i$. Choosing
$u_i=\alpha_i$ as basis of $\h^*$, so that the dual basis $\{u^i\}$ of $\h$ is
given by the fundamental coweights $\{\cow{i}\}$ yields $du_i=-dz_i/z_i$
and
$$du_i\,\tau(u^i)=-\frac{dz_i}{z_i}\,\tau(\cow{i})$$

\noindent Further, if $\alpha=\sum_i m_\alpha^i\alpha_i$ is a positive root,
then $e^\alpha=\prod_i z_i^{-m_\alpha^i}$ so that
$$\etaa{\alpha}=
\frac{e^{-\alpha}}{1-e^{-\alpha}}\,d\alpha=
-\sum_i m_\alpha^i
\frac{z_i^{m_\alpha^i-1}\prod_{j\neq i}z_j^{m_\alpha^j}}
{1-\prod_j z_j^{m_\alpha^j}}\,dz_i$$
which is a regular on each $T_i$.
It follows that $\res_{T_i}\Omega_1=0$ and
\begin{align*}
\res_{T_i}\Omega_2
&=
\imath_i^*\sum_\alpha\etaa{\alpha}\,[t_\alpha,\tau(\cow{i})]\\
\res_{T_i}\Omega_3
&=
\imath_i^*\sum_{j\neq i} \frac{dz_j}{z_j}\,[\tau(\cow{i}),\tau(\cow{j})]
\end{align*}
where $\imath_i$ is the inclusion $T_i\hookrightarrow\ol{T}$.

\subsection{}

Thus, for any $j\neq i$,
$$\res_{T_j\cap T_i}\res_{T_i}\Omega=[\tau(\cow{j}),\tau(\cow{i})]$$
which shows the necessity of the relations \eqref{eq:comm}.

\subsection{}\label{ss:end nec}

Let now $\alpha=\sum_i m_\alpha^i\alpha_i$ be a positive root and let $\ol
{T}_\alpha=\{\prod_iz_i^{m_\alpha^i}=1\}$ be the closure of $T_\alpha$
in $\ol{T}$. The intersection
$$T_{\alpha,i}=\ol{T}_\alpha\cap T_i\subset\ol{T}$$
is clearly nonempty if, and only if $\alpha(\cow{i})=m_\alpha^i=0$. When
that is the case, $T_{\alpha,i}$ is connected and contained in no other
$\ol{T}_\beta$, $\beta\neq j$ or $T_j$ for $j\neq i$.

It follows that whenever $\alpha(\cow{i})=0$,
$$\res_{T_{\alpha,i}}\res_{T_\alpha}\Omega=[t_\alpha,\tau(\cow{i})]$$
thus showing the necessity of \eqref{eq:cross} for $\alpha$ simple and
$w=1$. The general case follows by repeating the computations of the
last two subsections in the compactification of $T$ corresponding to a
different basis $w\Delta$ of simple roots and using the alternative form
of the connection $\nabla$ given by Proposition \ref{pr:chamber}.

\subsection{} 

We next turn to the sufficiency of the relations \eqref{eq:quadratic}--\eqref
{eq:comm}. This may be proved by embedding $T$ in the toric variety
corresponding to the fan determined by the chambers of $\Phi$ in $E$
and using a general integrability criterion of E. Looijenga as in \cite
[\S 1--2]{Lo}\footnote{note however that line 2 of the statement of Corollary
1.3 in \cite{Lo} should read "for every {\bf irreducible component $I$ of a}
codimension two intersection", the words in bold are missing in \cite{Lo}.}.
We prefer a more direct approach which will occupy \S \ref{ss:Omega 2}--
\ref{ss:end flat}.

\subsection{}\label{ss:Omega 2}

Since the relations \eqref{eq:comm} imply that $\Omega_3=0$,
we need to show that the relations \eqref{eq:quadratic}--\eqref{eq:cross}
imply that $\Omega_1+\Omega_2=0$. To this end, we rewrite $\Omega_2$
in a different form below and, in \S \ref{eq:Omega 1} rewrite $\Omega_1$.

\begin{lemma}\label{le:Omega 2}
Modulo the relations \eqref{eq:cross}, the curvature term
$$\Omega_2^\alpha=\sum_i \etaa{\alpha}\wedge du_i\,[t_\alpha,\tau(u^i)]$$
corresponding to $\alpha\in\Phi_+$ is equal to
$$\sum_{\beta\in\Phi_+\cap w\Phi_-}
\etaa{\alpha}\wedge d\beta\,[t_\alpha,t_\beta]$$
for any $w\in W$ such that $w^{-1}\alpha$ is a simple root.
\end{lemma}
\proof
Let $w\in W$ be such that $w^{-1}\alpha$ is a simple root $\alpha_i$.
By \eqref{eq:chamber},
$$\Omega_2^\alpha=
\sum_j \etaa{\alpha}\wedge du_j\,
[t_\alpha,\tau_w(u^j)+\sum_{\beta\in\Phi_+\cap w\Phi_-}\beta(u^j)\,t_\beta]$$
Choosing $u_j=w\alpha_j$ yields a commutator $[t_\alpha,\tau_w(u^j)]=
[t_\alpha,\tau_w(w\cow{j})]$ which is zero for all $j\neq i$ by \eqref{eq:cross}.
Since $d\alpha\wedge w\alpha_i=0$, this yields
$$\Omega_2^\alpha=
\sum_j \etaa{\alpha}\wedge du_j\,
[t_\alpha,\sum_{\beta\in\Phi_+\cap w\Phi_-}\beta(u^j)\,t_\beta]=
\sum_{\beta\in\Phi_+\cap w\Phi_-}\etaa{\alpha}\wedge d\beta\,
[t_\alpha,t_\beta]$$
since $\beta=\beta(u^i)u_i$ \halmos

\subsection{}\label{eq:Omega 1}

For any $a\in Q$, let $\eta_a$ be the meromorphic one--form on $T$
given by
$$\eta_a=\frac{da}{e^a-1}$$

\begin{lemma}\label{le:spillover}
Let $\Psi\subset\Phi$ be a rank 2 root subsystem and consider the
curvature term
$$\Omega_1^\Psi=
\half{1}\sum_{\alpha,\beta\in\Psi_+}
\eta_\alpha\wedge\eta_\beta\,[t_\alpha,t_\beta]$$
Assume that the relations \eqref{eq:quadratic} hold. Then, if $\Psi$ is of type
$\sfA_1\times\sfA_1$,
\begin{align}
\Omega_1^\Psi&=0
\intertext{If $\Psi$ is of type $\sfA_2$ with $\Psi_+=\{\alpha,\beta,\alpha+\beta\}$}
\Omega_1^\Psi&=
-\eta_{\alpha+\beta}\wedge d\beta\, [t_{\alpha+\beta},t_\beta]
\intertext{If $\Psi$ is of type $\sfB_2$ with $\Psi_+=\{\alpha,\beta,\alpha\pm\beta\}$}
\Omega_1^\Psi&=
-\eta_{\alpha+\beta}\wedge d\beta\, [t_{\alpha+\beta},t_\beta]
-\eta_\alpha\wedge d(\alpha-\beta)\, [t_{\alpha},t_{\alpha-\beta}]
\intertext{If $\Psi$ is of type $\sfG_2$, with $\Psi_+=\{\alpha_1,\alpha_2,\alpha_1+
\alpha_2,\alpha_1+2\alpha_2,\alpha_1+3\alpha_2,2\alpha_1+3\alpha_2\}$}
\Omega_1^\Psi
=&-
\eta_{\alpha_1+\alpha_2}\wedge d\alpha_1\,
[t_{\alpha_1+\alpha_2},t_{\alpha_1}]\label{eq:G2}\\
&-
\eta_{\alpha_1+2\alpha_2}\wedge\left(
d\alpha_2\,[t_{\alpha_1+2\alpha_2},t_{\alpha_2}]+
d(\alpha_1+3\alpha_2)\,[t_{\alpha_1+2\alpha_2},t_{\alpha_1+3\alpha_2}]\right)
\nonumber\\
&-
\eta_{\alpha_1+3\alpha_2}\wedge d\alpha_2\,
[t_{\alpha_1+3\alpha_2},t_{\alpha_2}]
\nonumber\\
&-
\eta_{2\alpha_1+3\alpha_2}\wedge\left(
d\alpha_1\,
[t_{2\alpha_1+3\alpha_2},t_{\alpha_1}]+
d(\alpha_1+\alpha_2)\,
[t_{2\alpha_1+3\alpha_2},t_{\alpha_1+\alpha_2}]\right)
\nonumber
\end{align}
\end{lemma}
\proof
If $\Psi$ is of type $\sfA_1\times\sfA_1$, the result follows from \eqref{eq:A1A1}.

For other types, we shall need the following easily verified identity. For $a,b\in
Q$, set
$$\eta_{a,b}=\frac{da\wedge db}{e^{a+b}-1}$$
Then,
\begin{equation}\label{eq:trigo OS}
\eta_a\wedge\eta_b=
\eta_a\wedge\eta_{a+b}+
\eta_{a+b}\wedge\eta_b+
\eta_{a,b}
\end{equation}
We shall apply \eqref{eq:trigo OS} to $\Omega_1^\Psi$ repeatedly, specifically to
terms of the form $\eta_\alpha\wedge\eta_\beta$ with $\alpha+\beta\in\Phi_+$,
until no such terms are left.

For $\Psi_+=\{\alpha,\beta,\alpha+\beta\}$ of type $\sfA_2$, this yields
$$\Omega_1^\Psi=
\eta_{\alpha+\beta}\wedge\eta_\alpha\,[t_{\alpha+\beta}+t_\beta,t_\alpha]+
\eta_{\alpha+\beta}\wedge\eta_\beta\,[t_{\alpha+\beta}+t_\alpha,t_\beta]
+\eta_{\alpha,\beta}\,[t_\alpha,t_\beta]
$$
By \eqref{eq:Phi}, the first two commutators are 0 and the third is
equal to $[t_\alpha,t_\beta]=[t_\alpha+t_\beta,t_\beta]=-[t_{\alpha+\beta},
t_\beta]$. This yields the required answer since
\begin{equation}\label{eq:eta eta}
\eta_{a,b}=
\eta_{a+b}\wedge db=-\eta_{a+b}\wedge da
\end{equation}

To keep track of the repeated applications of \eqref{eq:trigo OS}
for $\Psi$ ot type $\sfB_2,\sfG_2$, we proceed as follows. Recall
that the height of $\alpha=\sum_im_\alpha^i\alpha_i\in\Phi_+$ is
defined by $\height(\alpha)=\sum_im_\alpha^i$. Arrange pairs of distinct
roots $(\alpha,\beta)$ on consecutive rows according to the value
of $\height(\alpha)+\height(\beta)$: each $(\alpha,\beta)$ stands for a term
$\eta_\alpha\wedge\eta_\beta$. From each pair $(\alpha,\beta)$
such that $\alpha+\beta\in\Phi_+$ draw an arrow to $(\alpha,\alpha
+\beta)$ and $(\alpha+\beta,\beta)$ to signify that \eqref{eq:trigo OS}
has been applied with $a=\alpha$ and $b=\beta$.

For $\Psi_+=\{\alpha,\beta,\alpha\pm\beta\}$ of type $\sfB_2$ with
simple roots $\alpha_1=\alpha-\beta$, $\alpha_2=\beta$, the
corresponding graph is
{\small 
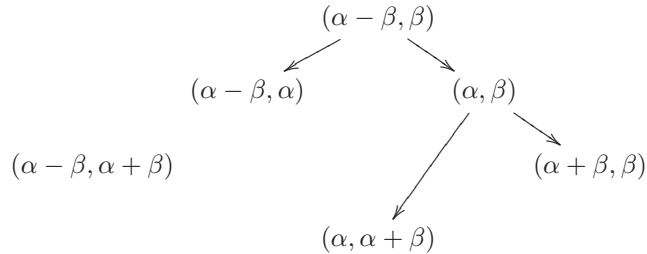
$$\begin{xy}
\xymatrix@=0pt@C=0pt@R=12pt{
			&&(\alpha-\beta,\beta)\ar[dl]\ar[dr]&& \\
&(\alpha-\beta,\alpha)	&		&(\alpha,\beta)\ar[dr]\ar[ddl]&\\
(\alpha-\beta,\alpha+\beta)&	&&	&(\alpha+\beta,\beta)\\
		&&(\alpha,\alpha+\beta)&&		}
\end{xy}$$}
This yields an "$\eta_a\wedge\eta_b$" component of $\Omega_1^\Psi$
equal to
\begin{multline*}
\eta_{\alpha-\beta}\wedge\eta_\alpha
\,[t_{\alpha-\beta},t_\alpha+t_\beta]\\
+\etaeta{\alpha-\beta}{\alpha+\beta}
+\eta_{\alpha+\beta}\wedge\eta_\beta
\,[t_{\alpha+\beta}+t_\alpha+t_{\alpha-\beta},t_\beta]\\
+
\eta_\alpha\wedge\eta_{\alpha+\beta}\,
\left([t_\alpha,t_{\alpha+\beta}+t_\beta]+[t_{\alpha-\beta},t_\beta]\right)
\end{multline*}
The second commutator is equal to zero by \eqref{eq:A1A1B2}, the first
by \eqref{eq:Phi} and \eqref{eq:A1A1B2} and the third by \eqref{eq:Phi}.
By \eqref{eq:Phi}, the coefficient of $\eta_\alpha\wedge\eta_{\alpha+\beta}$
is equal to $-[t_\alpha,t_{\alpha-\beta}]+[t_{\alpha-\beta},t_\beta]=[t_{\alpha-
\beta},t_\alpha+t_\beta]=0$.

$\Omega_1^\Psi$ is therefore equal to its "$\eta_{a,b}$" component, namely
$$\eta_{\alpha-\beta,\beta}\,[t_{\alpha-\beta},t_{\beta}]+
\eta_{\alpha,\beta}
\,[t_{\alpha}+t_{\alpha-\beta},t_{\beta}]$$
which yields the required answer since, by \eqref{eq:Phi}
$$[t_{\alpha}+t_{\alpha-\beta},t_{\beta}]=
-[t_{\alpha+\beta},t_{\beta}]$$
and by \eqref{eq:Phi} and \eqref{eq:A1A1B2}
$$[t_{\alpha-\beta},t_{\beta}]=
-[t_{\alpha-\beta},t_{\alpha}+t_{\alpha+\beta}]=
-[t_{\alpha-\beta},t_{\alpha}]$$
while, as previously noted
$$\eta_{a,b}=\eta_{a+b}\wedge db=-\eta_{a+b}\wedge da$$

Assume now that $\Psi$ is of type $\sfG_2$ and has simple roots $\alpha
_1, \alpha_2$, with $\alpha_1$ long. The sets of long and short positive
roots are, respectively
$$\Psi_+^\ell=\{\alpha_1,2\alpha_1+3\alpha_2,\alpha_1+3\alpha_2\}
\aand
\Psi_+^s=\{\alpha_1+\alpha_2,\alpha_1+2\alpha_2,\alpha_2\}$$
and the pairs $(\beta,\gamma)$ of orthogonal positive roots are
$$(\alpha_1,\alpha_1+2\alpha_2),\qquad
(2\alpha_1+3\alpha_2,\alpha_2),\qquad
(\alpha_1+3\alpha_2,\alpha_1+\alpha_2)
$$
The corresponding graph reads
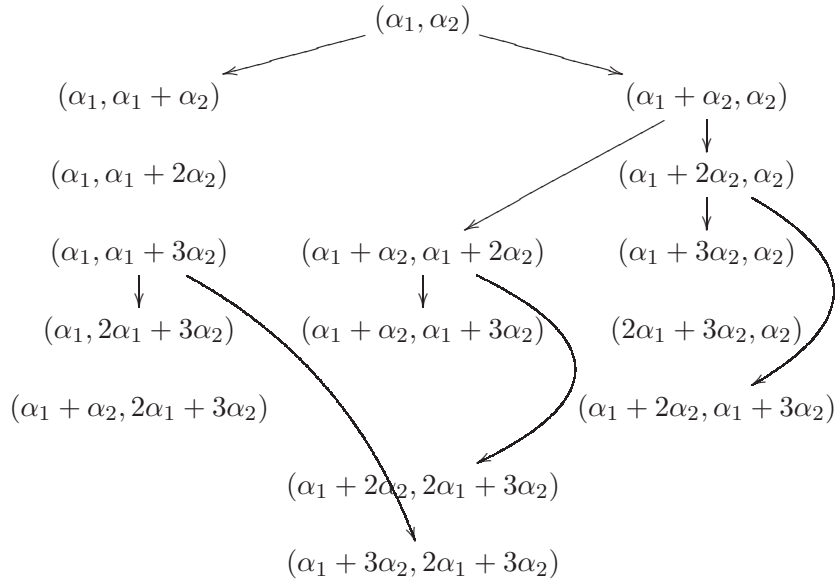
{$$
\xymatrix@=0pt@C=0pt@R=12pt{
&(\alpha_1,\alpha_2)\ar[dl]\ar[dr]&\\
(\alpha_1,\alpha_1+\alpha_2)&&(\alpha_1+\alpha_2,\alpha_2)\ar[ddl]\ar[d]\\
(\alpha_1,\alpha_1+2\alpha_2)&&(\alpha_1+2\alpha_2,\alpha_2)\ar[d]\ar@/^4pc/[ddd]\\
(\alpha_1,\alpha_1+3\alpha_2)\ar[d]\ar@/^1.5pc/[ddddr]&(\alpha_1+\alpha_2,\alpha_1+2\alpha_2)\ar[d] \ar@/^4.8pc/[ddd]&(\alpha_1+3\alpha_2,\alpha_2)\\
(\alpha_1,2\alpha_1+3\alpha_2)&(\alpha_1+\alpha_2,\alpha_1+3\alpha_2)&(2\alpha_1+3\alpha_2,\alpha_2)\\
(\alpha_1+\alpha_2,2\alpha_1+3\alpha_2)&&(\alpha_1+2\alpha_2,\alpha_1+3\alpha_2)\\
&(\alpha_1+2\alpha_2,2\alpha_1+3\alpha_2)&\\
&(\alpha_1+3\alpha_2,2\alpha_1+3\alpha_2)&
}$$}
This yields an $\eta_{a,b}$ component of $\Omega_1^\Psi$ equal to
\begin{align*}
&
\eta_{\alpha_1,\alpha_2}\,
[t_{\alpha_1},t_{\alpha_2}]\\
&+
\eta_{\alpha_1+\alpha_2,\alpha_2}\,
[t_{\alpha_1+\alpha_2}+t_{\alpha_1},t_{\alpha_2}]\\
&+
\eta_{\alpha_1+2\alpha_2,\alpha_2}\,
[t_{\alpha_1+2\alpha_2}+t_{\alpha_1+\alpha_2}+t_{\alpha_1},t_{\alpha_2}]\\
&+
\eta_{\alpha_1,\alpha_1+3\alpha_2}\,
[t_{\alpha_1},t_{\alpha_1+3\alpha_2}]\\
&+
\eta_{\alpha_1+\alpha_2,\alpha_1+2\alpha_2}\,
\left([t_{\alpha_1+\alpha_2},t_{\alpha_1+2\alpha_2}+
t_{\alpha_2}]+[t_{\alpha_1},t_{\alpha_2}]\right)
\end{align*}
By \eqref{eq:variant}, the first commutator is equal to $-[t_{\alpha_1},t_{\alpha_1
+\alpha_2}]$. By \eqref{eq:Phi} and \eqref{eq:A1A1G2}, the third commutator is
equal to $-[t_{\alpha_1+3\alpha_2},t_{\alpha_2}]$. The second commutator is
therefore equal to
$$-[t_{\alpha_1+3\alpha_2}+t_{\alpha_1+2\alpha_2},t_{\alpha_2}]=
-[t_{\alpha_1+2\alpha_2},t_{\alpha_2}+t_{\alpha_1+3\alpha_2}]$$
where we used \eqref{eq:variant}. By \eqref{eq:A2G2}, the fourth commutator
is equal to $-[t_{\alpha_1},t_{2\alpha_1+3\alpha_2}]$. Finally, the coefficient of
$\eta_{\alpha_1+\alpha_2,\alpha_1+2\alpha_2}$ is equal to
$$-[t_{\alpha_1+\alpha_2},
t_{\alpha_1}+t_{2\alpha_1+3\alpha_2}]+[t_{\alpha_1},t_{\alpha_2}]=
[t_{2\alpha_1+3\alpha_2},t_{\alpha_1+\alpha_2}]$$
where we used \eqref{eq:variant}. Using \eqref{eq:eta eta} shows that the
\rhs of \eqref{eq:G2} is equal to the $\eta_{a,b}$--component of $\Omega
_1^\Psi$. A similar use of relations \eqref{eq:Phi}, \eqref{eq:A2G2} and 
\eqref{eq:A1A1G2} shows that the $\eta_a\wedge\eta_b$ component of
$\Omega_1^\Psi$ is zero and therefore completes the proof \halmos

\subsection{} 

\begin{corollary}\label{co:Omega 1}
Assume that the relations \eqref{eq:quadratic} and \eqref{eq:cross}
hold. Then, for any rank 2 root subsystem $\Psi\subset\Phi$, the
curvature term $\Omega_1^\Psi$ is equal to
$$\Omega_1^\Psi=
-\sum_{\alpha\in\Psi_+}\sum_{\beta\in\Psi_+\cap w_\alpha\Psi_-}
\eta_\alpha\wedge d\beta\,[t_\alpha,t_\beta]$$
where $w_\alpha$ is any element of the Weyl group of $\Psi$ such
that $w_\alpha^{-1}\alpha$ is simple in $\Psi$.
\end{corollary}
\proof Lemma \ref{le:spillover} and a simple \cbc inspection show that
$\Omega_1^\Psi$ does indeed have the above form for well--chosen
elements $w_\alpha$ (specifically, $w_\alpha$ should be an element
of shortest length such that $w_\alpha^{-1}\alpha$ is simple in $\Psi$).
But Lemma \ref{le:Omega 2} applied to $\Psi$ implies that the expression
$$\sum_{\beta\in\Psi_+\cap w_\alpha\Phi_-}
\eta_\alpha\wedge d\beta\,[t_\alpha,t_\beta]$$
is independent of the choice of $w_\alpha$, hence the conclusion. \halmos

\subsection{Completion of the proof of (1) of Theorem \ref{th:trigo flat}}\label{ss:end flat}

Fix $\alpha\in\Phi_+$ and let $\R_2(\alpha)$ be the set of complete,
rank 2 subsystems $\Psi\subset\Phi$ containing $\alpha$ as a non--simple
root. For $\Psi\in \R_2(\alpha)$, denote the corresponding Weyl group
by $W(\Psi)$.

For any $w\in W$, denote by $N(w)\subset\Phi_+$ the set
$$N(w)=\{\beta\in\Phi_+|\, w\beta\in\Phi_-\}$$
and define similarly $N_\Psi(w)\subset\Psi_+$ for any $w\in W(\Psi)$.
 
The following result, together with Lemma \ref{le:Omega 2} and Corollary
\ref{co:Omega 1} show that $\Omega_1+\Omega_2=0$ and therefore
that the connection \eqref{eq:trigo} is flat.
 
\begin{proposition}
Let $w\in W$ be such that $w^{-1}\alpha$ is simple in $\Phi$.
\begin{enumerate}
\item For any $\Psi\in\R_2(\alpha)$, the intersection $N(w^{-1})\cap\Psi$
is non--empty.
\item The following holds
$$N(w^{-1})=\bigsqcup_{\Psi\in\R_2(\alpha)}N(w^{-1})\cap\Psi$$
\item For any $\Psi\in \R_2(\alpha)$, there exists a unique $w_\Psi\in W
(\Psi)$ such that
$$N(w^{-1})\cap\Psi=N_\Psi(w_\Psi^{-1})$$
\end{enumerate}

\end{proposition}
\proof For any pair of non--proportional roots $\beta,\gamma\in\Phi$, let
$\<\beta,\gamma\>\subset\Phi$ be the complete, rank 2 subsystem generated
by $\beta$ and $\gamma$. We claim that the map $\beta\to\<\alpha,\beta\>$
induces a bijection
$$\rho:\left.N(w^{-1})\right/_\sim\longrightarrow \R_2(\alpha)$$
where $\sim$ is the equivalence relation defined by $\beta\sim\beta'$
if $\<\alpha,\beta\>=\<\alpha,\beta'\>$. This clearly proves (1) and (2)
since $\rho^{-1}(\Psi)=\Psi\cap N(w^{-1})$.

To see this, we shall need some notation. For any complete subsystem
$\Psi\subset\Phi$, let $\Psi^\perp=\bigcap_{\beta\in\Psi}\H_\beta\subset
E$, where $\H_\beta=\Ker(\beta)$. Set $E_\Psi=E/\Psi^\perp$ so that $
\Psi$ may be regarded as a root system in $E_\Psi^*$, and let $\pi_\Psi:
E\to E_\Psi$ be the corresponding projection. If $\beta\in N(w^{-1})$, the
wall $\H_\beta$ separates the fundamental chamber $C\subset E$ of $
\Phi$ and $C'=w(C)$. Thus, if $\Psi=\<\alpha,\beta\>$, $\pi_\Psi(\H_\beta)$
separates the fundamental chamber $\pi_\Psi(C)$ of $\Psi$ and $\pi_\Psi
(C')$. It follows that $\pi_\Psi(C)\neq\pi_\Psi(C')$ so that $\alpha$ is not
simple in $\Psi$ since $\pi_\Psi(\H_\alpha)$ is a wall of $\pi_\Psi(C')$,
whence $\<\alpha,\beta\>\in \R_2(\alpha)$ and $\rho$ is a well--defined
embedding.

It is also surjective since if $\Psi\in \R_2(\alpha)$, there exists a $\beta\in
\Psi_+$ such that $\pi_\Psi(\H_\beta)$ separates $\pi_\Psi(C)$ and $\pi_
\Psi(C')$ so that $\H_\beta$ separates $C$ and $C'$ and therefore lies in
$N(w^{-1})$.

Finally, for a given $\Psi\in \R_2(\alpha)$, the set $N(w^{-1})\cap\Psi$
consists of those $\beta\in\Psi_+$ which separate $C$ and $C'$ and
therefore $\pi_\Psi(C)$ and $\pi_\Psi(C')$. It is therefore equal to
$\Psi_+\cap w_\Psi\Psi_-$ where $w_\Psi\in W(\Psi)$ is the unique
element such that $\pi_\Psi(C')=w_\Psi\pi_\Psi(C)$. \halmos

\subsection{}

We now turn to part (2) of Theorem \ref{th:trigo flat}. We shall need the
following

\begin{lemma}\label{le:alternative}
The relations \eqref{eq:quadratic} imply that the following holds for any
$\alpha\in\Phi_+$, $w\in W$ such that $w^{-1}\alpha$ is simple and $v
\in \h$ such that $\alpha(v)=0$
\begin{equation}\label{eq:twt}
[t_\alpha,\sum_{\beta\in\Phi_+}\sign(w^{-1}\beta)\beta(v)t_\beta]=0
\tag{$t ^wt$}
\end{equation}
where $\sign(\gamma)=\pm 1$ depending on whether $\gamma\in\pm\Phi_+$.
\end{lemma}
\proof
Let $\T$ be the algebra generated by symbols $t_\alpha$, $\alpha\in\Phi$
subject to the relations $t_{-\alpha}=t_\alpha$ and \eqref{eq:quadratic}. 
The Weyl group acts on $\T$ by $wt_\alpha=t_{w\alpha}$ and it is easy
to check that the above relation holds for a triple $(\alpha,w,v)$ if, and
only if, it holds for $(w^{-1}\alpha,1,w^{-1}v)$. We may therefore assume
that $\alpha$ is simple and that $w=1$. Since the \lhs of \eqref{eq:twt}
may be written as
$$\sum_{\Psi}\sum_{\beta\in\Psi_+}[t_\alpha,\beta(v)t_\beta]$$
where $\Psi$ ranges over the complete, rank 2 subsystems of $\Phi$
containing $\alpha$, it is sufficient to prove \eqref{eq:twt} when $\Phi$ is of
rank 2. In this case, it follows by a simple \cbc verification.\Omit{Indeed, if
$\Phi_+=\{\alpha_1,\alpha_2,\alpha_1+\alpha_2\}$ if of type $\sfA_2$,
then
$$[t_{\alpha_1},\sum_{\beta\in\Phi_+}\beta(\cow{2})t_\beta]=
[t_{\alpha_1},t_{\alpha_2}+t_{\alpha_1+\alpha_2}]=0$$
by \eqref{eq:Phi} and, by symmetry, $[t_{\alpha_2},\sum_\beta\beta
(\cow{1})t_{\beta}]=0$.

If $\Phi_+=\{\alpha_1,\alpha_2,\alpha_1+\alpha_2,\alpha_1+2\alpha_2\}$
if of type $\sfB_2$, then $\alpha=\alpha_1$ and $v=\cow{2}$ yield
$$[t_{\alpha_1},\sum_{\beta\in\Phi_+}t_\beta]+
[t_{\alpha_1},t_{\alpha_1+2\alpha_2}]=0$$
by \eqref{eq:Phi} and \eqref{eq:A1A1B2}, while $\alpha=\alpha_2$ and
$v=\cow{1}$ give $[t_{\alpha_2},\sum_{\beta\in\Phi_+}t_\beta]=0$ by
\eqref{eq:Phi}.
Finally,}
For example, if $\Phi_+=\{\alpha_1,\alpha_2,\alpha_1+\alpha_2,\alpha_1+2\alpha_2,
\alpha_1+3\alpha_2,2\alpha_1+3\alpha_2\}$ is of type $\sfG_2$, we have
$$[t_{\alpha_2},\sum_{\beta\in\Phi_+}\beta(\cow{1})t_\beta]=
[t_{\alpha_2},\sum_{\beta\in\Phi_+}t_\beta]+[t_{\alpha_2},t_{2\alpha_1+3\alpha_2}]=0$$
by \eqref{eq:Phi} and \eqref{eq:A1A1G2}, while
$$[t_{\alpha_1},\sum_{\beta\in\Phi_+}\beta(\cow{2})t_\beta]=
[t_{\alpha_1},\sum_{\beta\in\Phi_+}t_\beta]+
[t_{\alpha_1},t_{\alpha_1+2\alpha_2}]+
2[t_{\alpha_1},t_{2\alpha_1+3\alpha_2}+t_{\alpha_1+3\alpha_2}]$$
which is equal to zero by \eqref{eq:Phi}, \eqref{eq:A1A1G2} and \eqref{eq:A2G2}
\halmos

\subsection{}\label{ss:alternative}

Recall from \S \ref{ss:delta} that $\delta:\h\to A$ is defined by
$$\delta(v)=\tau(v)-\half{1}\sum_{\alpha\in\Phi_+}\alpha(v)t_\alpha$$
The following proves part (2) of Theorem \ref{th:trigo flat}

\begin{proposition}
Modulo the relations \eqref{eq:quadratic}, the relations \eqref
{eq:cross} are equivalent to
\begin{equation}
[t_\alpha,\delta(v)]=0
\tag{$t\delta$}
\end{equation}
for any $\alpha\in\Phi$ and $v\in \h$ such that $\alpha(v)=0$.
\end{proposition}
\proof
For any $w\in W$ \eqref{eq:chamber}, yields
$$\tau_w(v)=
\tau(v)-\negthickspace\negthickspace\negthickspace\negthickspace
\sum_{\alpha\in\Phi_+\cap w\Phi_-}\alpha(v)t_\alpha=
\delta(v)+
\half{1}\sum_{\beta\in\Phi_+}\sign(w^{-1}\beta)\beta(v)t_\beta(v)$$
The result now follows from Lemma \ref{le:alternative}. \halmos

\subsection{Equivariance under $W$}

Assume now that the algebra $A$ is acted upon by the Weyl group $W$
of $\Phi$.

\begin{proposition}\label{pr:equivariance}\hfill
\begin{enumerate}
\item The connection $\nabla$ is $W$--equivariant if, and only if
\begin{gather}
s_i(t_\alpha)=t_{s_i\alpha}\label{eq:equiv1}\\
s_i(\tau(x))-\tau(s_i x)=(\alpha_i,x)t_{\alpha_i}
\label{eq:equiv2}
\end{gather}
for any $\alpha\in\Phi$, simple reflection $s_i\in W$ and $x\in \h$.
\item Modulo \eqref{eq:equiv1}, the relation \eqref{eq:equiv2} is
equivalent to the $W$--equivariance of the linear map $\delta:\h
\to A$ defined by \eqref{eq:delta}.
\end{enumerate}
\end{proposition}

\proof (1) Since $s_i$ permutes the set $\Phi_+\setminus\{\alpha_i\}$
and, by \eqref{eq:trick}
$$\frac{1}{1-e^{-\alpha_i}}
=\frac{1}{e^{\alpha_i}-1}+1$$
we get
$$s_i^*\nabla=d-
\sum_{\alpha\in\Phi_+}\frac{d\alpha}{e^\alpha-1}s_i(t_\alpha)
-d\alpha_i\,s_i(t_{\alpha_i})
-s_i(\tau(u^j))d(s_i u_j)$$
Requiring that $s_i^*\nabla=\nabla$ and taking residues along each
subtorus $\{e^\alpha=1\}$ yields \eqref{eq:equiv1}. To proceed, note
that
$$s_i(\tau(u^j))d(s_i u_j)=s_i(\tau(s_i u^j))d(u_j)$$
since $\tau(u^j)du_j$ is independent of the choice of the dual bases $\{u^j\},
\{u_j\}$. Thus, $s_i^*\nabla=\nabla$ reduces to
$$\tau(u^j)du_j=s_i(\tau(s_i u^j))du_j+t_{\alpha_i}d\alpha_i$$
which, upon being contracted along the tangent vector $u^k$
yields \eqref{eq:equiv2} with $x=s_iu^j$. 

(2) It is easy to check that the map $\ol{\tau}(x)=1/2\sum_{\alpha\in\Phi
_+}(x,\alpha)t_\alpha$ satisfies \eqref{eq:equiv2}. The result now follows
since any two maps $\tau_i:\h\to A$ satisfying \eqref{eq:equiv2} differ by
a $W$--equivariant map. \halmos

\subsection{Flatness and equivariance}

The following is a direct corollary of Theorem \ref{th:trigo flat} and Proposition
\ref{pr:equivariance}

\begin{theorem}\label{th:W trigo flat}
The trigonometric connection
$$\nabla=
d-\sum_{\alpha\in\Phi_+}\frac{d\alpha}{e^\alpha-1}\,t_\alpha-d(u_i)\,\tau(u^i)$$
is flat and $W$--equivariant if, and only if the following relations hold
\vskip .1cm
\begin{itemize}
\item For any rank 2 root subsystem $\Psi\subset\Phi$ and $\alpha\in\Psi$,
$$[t_\alpha,\sum_{\beta\in\Psi_+}t_\beta]=0$$
\item For any $u,v\in \h$,
$$[\tau(u),\tau(v)]=0$$
\item For any simple root $\alpha_i$ and $u\in\Ker(\alpha_i)$,
$$[t_{\alpha_i},\tau(u)]=0$$
\item For any $\alpha\in\Phi$ and simple reflection $s_i\in W$,
$$s_i(t_\alpha)=t_{s_i\alpha}$$
\item For any $\alpha\in\Phi$ and $u\in\h$,
$$s_i(\tau(u))-\tau(s_i u)=(\alpha_i,u)t_{\alpha_i}$$
\end{itemize}
\end{theorem}

\begin{remark}
Theorem \ref{th:W trigo flat} was first proved by Cherednik in the special
case when $t_\alpha$ is equal to the orthogonal reflection $s_\alpha
\in W$ and shown to lead to the definition of the \daha of $W$ \cite
{Ch1,Ch2}.
\end{remark}

\section{The trigonometric Casimir connection}\label{se:trigo g}

\subsection{The Yangian $Y(\g)$ \cite{Dr1}}\label{ss:J(x)}

Let $\g$ be a finite--dimensional, simple Lie algebra over $\IC$ and $(\cdot,
\cdot)$ a non--degenerate, invariant bilinear form on $\g$. Let $\hbar$ be a
formal variable.
The Yangian $Y(\g)$ is the associative algebra over $\IC[\hbar]$ generated
by elements $x,J(x)$, $x\in\g$ subject to the relations
\begin{gather*}
\lambda x+\mu y\quad\text{(in $Y(\g)$) }=\lambda x+\mu y\quad\text{(in $\g$)}\\[1.2ex]
xy-yx=[x,y]\\[1.2ex]
J(\lambda x+\mu y)=\lambda J(x)+\mu J(y)\\[1.2ex]
[x,J(y)]=J([x,y])\\[1.2ex]
[J(x),J([y,z])]+[J(z),J([x,y])]+[J(y),J([z,x])]=\phantom{123456789}\phantom{123456789}\\
\phantom{123456789}\phantom{123456789}\hbar^2
([x,x_a],[[y,x_b],[z,x_c]])
\{x^a,x^b,x^c\}\\[1.2ex]
[[J(x),J(y)],[z,J(w)]]+[[J(z),J(w)],[x,J(y)]]=\phantom{123456789}\phantom{123456789}\\
\phantom{123456789}\phantom{123456789}\hbar^2
([x,x_a],[[y,x_b],[[z,w],x_c]])
\{x^a,x^b,J(x^c)\}
\end{gather*}
for any $x,y,z,w\in\g$ and $\lambda,\mu\in\IC$, where $\{x_a\},\{x^a\}$
are dual bases of $\g$ \wrt $(\cdot,\cdot)$ and
$$\{z_1,z_2,z_3\}=
\frac{1}{24}\sum_{\sigma\in\mathfrak{S}_3}
z_{\sigma(1)}z_{\sigma(2)}z_{\sigma(3)}$$
$Y(\g)$ is an $\IN$--graded $\IC[\hbar]$--algebra provided one sets $\deg(x)=0$,
$\deg(J(x))=1$ and $\deg(\hbar)=1$.

\subsection{Drinfeld's new realisation of $Y(\g)$ \cite{Dr2}}\label{ss:new realisation}

Let $\h\subset\g$ be a Cartan subalgebra of $\g$, $\Phi=\{\alpha\}\subset
\h^*$ the corresponding root system. Let $\{\alpha_i\}_{i\in\bfI}$ be a basis
of simple roots of $\Phi$ and $a_{ij}=2(\alpha_i,\alpha_j)/(\alpha_i,\alpha_i)$
the entries of the Cartan matrix $\bfA$ of $\g$. Set $d_i=(\alpha_i,\alpha_i)/2$,
so that $d_ia_{ij}=d_j a_{ji}$ for any $i,j\in\bfI$.

Let $\nu:\h\to\h^*$ be the isomorphism determined by the inner product
$(\cdot,\cdot)$ and set $t_i=\nu^{-1}(\alpha_i)=d_i\alpha_i^\vee$. For any
$i\in\bfI$, choose root vectors $x^\pm_i\in\g_{\pm\alpha_i}$ such that $[x_i
^+,x_i^-]=t_i$. Recall that $\g$ has a (slightly non--standard) presentation
in terms of the generators $t_i,x_i^\pm$ with relations
\begin{gather*}
[t_i,t_j]=0\\
[t_i,x_j^\pm]=\pm d_ia_{ij}x_j^\pm\\
[x_i^+,x_j^-]=\delta_{ij}t_i\\
\ad(x_i^\pm)^{1-a_{ij}}x_j^\pm=0
\end{gather*}

The Yangian $Y(\bfA)$ is the associative algebra over $\IC[\hbar]$ with
generators $X_{i,r}^\pm,T_{i,r}$, $i\in\bfI$, $r\in\IN$ and relations
\begin{gather*}
[T_{i,r},T_{j,s}]=0\\[1.1ex]
[T_{i,0},X_{j,s}^\pm]=\pm d_i a_{ij}X_{j,s}^\pm\\[1.1ex]
[T_{i,r+1},X_{j,s}^\pm]-[T_{i,r},X_{j,s+1}^\pm]=
\pm\half{\hbar}d_i a_{ij}(T_{i,r}X_{j,s}^\pm+X_{j,s}^\pm T_{i,r})\\[1.1ex]
[X_{i,r}^+,X_{j,s}^-]=\delta_{ij}T_{i,r+s}\\[1.1ex]
[X_{i,r+1}^\pm,X_{j,s}^\pm]-[X_{i,r}^\pm,X_{j,s+1}^\pm]=
\pm\half{\hbar}d_i a_{ij}(X_{i,r}^\pm X_{j,s}^\pm+X_{j,s}^\pm X_{i,r}^\pm)\\[1.1ex]
\sum_\pi
[X_{i,r_{\pi(1)}}^\pm,[X_{i,r_{\pi(2)}}^\pm,[\cdots,[X_{i,r_{\pi(m)}}^\pm,X_{j,s}^\pm]\cdots]]=0
\end{gather*}
where $i\neq j$ in the last relation, $m=1-a_{ij}$, $r_1,\ldots,r_m\in\IN$ is any
sequence of non--negative integers, and the sum is over all permutations $\pi$
of $\{1,\ldots,m\}$. $Y(\bfA)$ is $\IN$--graded by $\deg(T_{i,r})=r=\deg(X^\pm
_{i,r})$ and $\deg(\hbar)=1$.

\subsection{Isomorphism between the two presentations \cite{Dr2}}\label{ss:old new}

Choose root vectors $x_\alpha\in\g_\alpha$ for any $\alpha\in\Phi$
such that $(x_\alpha,x_{-\alpha})=1$ and let 
\begin{equation}\label{eq:truncated}
\kappa_\alpha=x_\alpha x_{-\alpha}+x_{-\alpha}x_\alpha
\end{equation}
be the truncated Casimir operator of the $\sl{2}$--subalgebra of
$\g$ corresponding to $\alpha$. Then, the assignment
\begin{gather*}
\varphi(t_i)=T_{i,0},\quad
\varphi(x^\pm_i)=X_{i,0}^\pm\\[1.2 ex]
\varphi(J(t_i))=T_{i,1}+\hbar\varphi(v_i)\\[1.2ex]
\varphi(J(x^\pm_i))=X_i^\pm+\hbar\varphi(w_i^\pm)
\end{gather*}
extends to an isomorphism $\varphi:Y(\g)\to Y(\bfA)$, where
\begin{align*}
v_i
&=
\frac{1}{4}\sum_{\beta\in\Phi_+} (\alpha_i,\beta)
\kappa_\alpha
-t_i^2/2\\[1.2ex]
w_i^\pm
&=
\pm\frac{1}{4}\sum_{\beta\in\Phi_+}
\left([x_i^\pm,x_{\pm\beta}]x_{\mp\beta}+x_{\mp\beta}[x_i^\pm,x_{\pm\beta}]\right)
-\frac{1}{4}(x_i^\pm t_i+t_ix_i^\pm)
\end{align*}

\subsection{The $W$--equivariant embedding $\delta_{a,b}:\h\to Y(\g)$}

The Yangian $Y(\g)$ is acted upon by the Lie algebra spanned by the
elements $x\in\g$ and is an integrable $\g$--module under this action.
In particular, the zero--weight subalgebra $Y(\g)^\h$ is acted upon by
the Weyl group $W$ of $\g$. Moreover, for any $a,b\in\IC$, the linear
map
\begin{equation}\label{eq:embedding}
\delta_{a,b}:\h\to Y(\g)^\h,\qquad
\delta_{a,b}(t)=at+bJ(t)
\end{equation}
is $W$--equivariant.

In terms of the new realisation of $Y(\g)$, the map $\delta_{a,b}$ becomes
\begin{gather}
\wt{\delta}_{a,b}=\varphi\circ\delta_{a,b}:\h\to Y(\bfA)^\h
\nonumber\\
\wt{\delta}_{a,b}(t)=
at+b\left(T(t)_1+\frac{\hbar}{4}\sum_{\beta\in\Phi_+}
\beta(t)\kappa_\beta-\frac{\hbar}{2}\sum_i(t,t^i)t_i^2\right)
\label{eq:2nd embedding}
\end{gather}
where $\{t^i=\cow{i}\}$ is the basis of $\h$ dual to $\{t_i\}$ given by
the fundamental coweights, $T(-)_1:\h\to Y(\bfA)$ is the embedding
$t\to (t,t^i)T_{i,1}$ and, deviating slightly from the notation of \S \ref
{ss:new realisation}, we identify $\g\subset Y(\g)$ with the Lie subalgebra
of $Y(\bfA)$ spanned by $T_{i,0},X_{i,0}^\pm$.

\subsection{The linear map $\tau_{a,b}:\h\to Y(\g)$}

Let $\epsilon\in\IC$. For any root $\alpha$, set
\begin{equation}
K_\alpha=\kappa_\alpha+\epsilon q_\alpha
\qquad\text{where}\qquad
q_\alpha=\frac{\nu^{-1}(\alpha)}{(\alpha,\alpha)}^2
\label{eq:K alpha}
\end{equation}
and $\kappa_\alpha$ is the truncated Casimir given by \eqref{eq:truncated}.
For any $a,b\in\IC$, define a map $\tau_{a,b}:\h\to Y(\g)^\h$ by
$$\tau_{a,b}(t)=
\half{\hbar}\sum_{\alpha\in\Phi_+}(t,\alpha)K_\alpha+\delta_{a,b}(t)$$
where $\delta_{a,b}$ is given by \eqref{eq:embedding}.

\begin{proposition}\label{pr:tau ab}\hfill
\begin{enumerate}
\item The elements $K_\alpha$ satisfy
$$K_{-\alpha}=K_\alpha\aand w(K_\alpha)=K_{w\alpha}$$
for any $w\in W$.
\item The following holds for any $t\in\h$ 
$$s_i(\tau_{a,b}(t))-\tau_{a,b}(s_i t)=\hbar(\alpha_i,t) K_{\alpha_i}$$
\item If $b=-2$, then $\wt{\tau}_{a,b}=\varphi\circ\tau_{a,b}:\h\to
Y(\bfA)^\h$ is given by
$$\wt{\tau}_{a,-2}(t)
=
at-2\left(T(t)_1-\half{\hbar}\sum_i(t,t^i)t_i^2\right)
+\epsilon\frac{\hbar}{2}\sum_{\alpha\in\Phi_+}\alpha(t)q_\alpha
$$
and satisfies in addition
$$[\wt{\tau}_{a,-2}(t),\wt{\tau}_{a,-2}(t')]=0$$
for any $t,t'\in\h$.
\item For any $\alpha\in\Phi_+$ and $t\in\h$ such that $\alpha(t)=0$
$$[K_\alpha,\delta_{a,b}(t)]=0$$
\end{enumerate}
\end{proposition}
\proof (1) is obvious. (2) follows from the second part of Proposition \ref
{pr:equivariance}, (1) and the $W$--equivariance of $\delta_{a,b}$. For
(3), we have by \eqref{eq:2nd embedding}
\begin{equation*}
\begin{split}
\wt{\tau}_{a,b}(t)
&=
\half{\hbar}\sum_{\alpha\in\Phi_+}\alpha(t)K_\alpha+
at+b\left(T(t)_1+\frac{\hbar}{4}\sum_{\beta\in\Phi_+}\beta(t)\kappa_\beta-\frac{\hbar}{2}\sum_i(t,t^i)t_i^2\right)\\
&=
at+
\frac{\hbar}{2}\sum_{\alpha\in\Phi_+}\alpha(t)(K_\alpha+\frac{b}{2}\kappa_\alpha)+
b(T(t)_1-\frac{\hbar}{2}\sum_i(t,t^i)t_i^2)
\end{split}
\end{equation*}
which, for $b=-2$ reduces to the claimed expression. The commutativity of $\wt{\tau}
_{a,-2}(t)$ and $\wt{\tau}_{a,-2}(t')$ then follows from that of of the $T_{i,1}$.
(4) follows easily from the defining relations of $Y(\g)$. \halmos

\subsection{}\label{ss:J & T}

For simplicity, we henceforth set $\epsilon=0=a$ and $b=-2$ in equations
\eqref{eq:embedding} and \eqref{eq:K alpha}, although Theorem \ref{th:trigo Casimir}
below is true for any values of $a,\epsilon$.
Thus, $K_\alpha=\kappa_\alpha$, $\delta=\delta_{0,-2}:\h\to Y(\g)$ is given by
$\delta(t)=-2J(t)$ and the corresponding map $\wt{\tau}=\wt{\tau}_{0,-2}:\h\to
Y(\bfA)^\h$ by
$$\wt{\tau}(t)=-2T(t)_1+\hbar (t,t^i)t_i^2$$

\subsection{The trigonometric Casimir connection of $\g$}

Let $H\reg\subset H$ be given by \eqref{eq:Hreg} and $\Y^\h$ the trivial
bundle over $H\reg$ with fibre $Y(\g)^\h\cong Y(\bfA)^\h$.

\begin{definition} The {\it trigonometric Casimir connection} of $\g$ is the
connection $\nabla$ on $\Y^\h$ given by either of the following forms
\begin{align}
\nabla&=
d-\half{\hbar}\sum_{\alpha\in\Phi_+}\frac{e^\alpha+1}{e^\alpha-1}d\alpha\,
\kappa_\alpha
+2du_i\,J(u^i)
\label{eq:J delta form}\\[1.5 ex]
&=
d-\hbar\sum_{\alpha\in\Phi_+}\frac{d\alpha}{e^\alpha-1}\,\kappa_\alpha
+2du_i\,\left(T(u^i)_1-\half{\hbar}(u^i,t^j)t_j^2\right)
\label{eq:T tau form}
\end{align}
whose equality follows from \S \ref{ss:delta} and \S \ref{ss:J & T}.
\end{definition}

\subsection{}

\begin{theorem}\label{th:trigo Casimir}
The trigonometric Casimir connection is flat and $W$--equivariant.
\end{theorem}
\proof By Theorem \ref{th:trigo flat}, Proposition \ref{pr:equivariance}
and Proposition \ref{pr:tau ab}, we need only check that the elements
$t_\alpha=\kappa_\alpha$ satisfy the relations \eqref{eq:quadratic}.
This follows as in \cite[Thm. 2.3]{MTL} and \cite[Thm 2.2]{TL1}.
Specifically, if $\Psi\subset\Phi$ is a rank 2 root subsystem, the sum
$\sum_{\alpha\in\Phi_+}\kappa_\alpha$ is, up to Cartan terms, the
Casimir operator for the rank 2 subalgebra $\g_\Psi\subset\g$
determined by $\Psi$ and therefore commutes with each summand
$\kappa_\alpha$, $\alpha\in\Phi_+$. \halmos
\newpage

\begin{remark}\label{rk:scaling} Since the relations of Theorem \ref
{th:trigo flat} and Proposition \ref{pr:equivariance} are homogeneous,
Theorem \ref{th:trigo Casimir} proves in fact the flatness and $W
$--equivariance of the one--parameter family of connections
\begin{align*}
\nabla&=
d-\lambda^{-1}\left(\half{\hbar}\sum_{\alpha\in\Phi_+}
\frac{e^\alpha+1}{e^\alpha-1}d\alpha\,
\kappa_\alpha
-2du_i\,J(u^i)\right)\\[1.5 ex]
&=
d-\lambda^{-1}\left(\hbar\sum_{\alpha\in\Phi_+}
\frac{d\alpha}{e^\alpha-1}\,\kappa_\alpha
-2du_i\,\left(T(u^i)_1-\half{\hbar}(u^i,t^j)t_j^2\right)\right)
\end{align*}
where $\lambda$ varies in $\IC^\times$.
\end{remark}

\section{The monodromy conjecture}\label{se:monodromy}

We show in this section that the monodromy of the trigonometric Casimir
connection $\nabla$ gives rise to representations of the affine braid group
$\wh{B}$ corresponding to $\g$ and give a conjectural description of it in
terms of the quantum Weyl group operators of the quantum loop algebra
$U_\hbar(L\g)$.

\subsection{Monodromy representation}

Since $H$ is of simply--connected type, the Weyl group $W$
acts freely on $H\reg$ and the fundamental group of the quotient $H\reg
/W$ is isomorphic to the affine braid group $\wh{B}$ \cite{NVD,vdL}.

Let $V$ be a \fd $Y(\g)$--module and $\IV$ the holomorphically trivial
vector bundle over $H\reg$ with fibre $V$. The connection $\nabla$
induces a flat connection on $\IV$. To push it down to the quotient by
$W$ we use the 'up and down' trick of \cite[p. 224]{MTL} to circumvent
the fact that $W$ does not in general act on $V$. To this end, we shall
need a few basic results about Tits extensions of (affine) Weyl groups
which are gathered in \S \ref{se:Tits}.

Specifically, since $V$ is an integrable $\g$--module, the triple
exponentials
$$\exp(e_{\alpha_i})\exp(-f_{\alpha_i})\exp(e_{\alpha_i})$$
defined by a choice of simple root vectors $e_{\alpha_i}\in\g_{\alpha_i},
f_{\alpha_i}\in\g_{-\alpha_i}$ are well--defined elements of $GL(V)$.
They give rise to an action on $V$ of an extension of $W$ by the sign
group $\IZ_2^{\dim\h}$ called the {\it Tits extension} $\wt{W}$ of $W$
(Definition \ref{de:Tits extension} and Prop. \ref{pr:Tits action}). By
Theorem \ref{th:reduced affine Tits},
$\wt{W}$ is a quotient of the affine braid group $\wh{B}$ which may
therefore be made to act on $V$. It is then easy to check that the
pull--back of the flat vector bundle $(\IV,\nabla)$ to the universal
cover of $H\reg$ is equivariant under $\wh{B}$ acting by deck
transformations on the base and through the $\wt{W}$--action
on the fibres.

\subsection{}

Let $L\g=\g[t,t^{-1}]$ be the loop algebra of $\g$ and $\Uhghat$
the corresponding quantum loop algebra, viewed as a topological
Hopf algebra over the ring of formal power series $\ICh$. Thus,
$\Uhghat$ has Chevalley generators $E_i,F_i$, where $i$ ranges
over the set $\wh{\bfI}=\bfI\sqcup\{0\}$ of nodes of the affine Dynkin
diagram of $\g$ and a Cartan subalgebra isomorphic to $\h$ and
spanned by $H_i$, $i\in\bfI$ and $H_0=-H_{\theta}=-\sum_{i\in\bfI}a_i
H_i$, where $\theta\in\h^*$ is the highest root and the integers
$a_i$ are given by $\theta^\vee=\sum_i a_i\root{i}^\vee$.

\subsection{}

By a {\it \fd representation} of $\Uhghat$ we shall mean a module
$\V$ which is topologically free and finitely--generated over $\IC\fml$.
Such a $\V$ is integrable and therefore endowed with a \qW action
of the affine braid group $\wh{B}$ \cite{Lu,KR,So}. This action is
given by letting the generator corresponding to $i\in\wh{\bfI}$ act
by
$$\ol{S}\ih\medspace v=
\sum_{\substack{a,b,c\in\IZ : \\a-b+c=-\lambda(\cor{i})}}
(-1)^{b}q_i^{b-ac}
E_i^{(a)}F_i^{(b)}E_i^{(c)}
v$$
where $v\in\V$ if of weight $\lambda\in\h^*$ and $X_i^{(a)}$ is the
divided power $X^a/[a]_i!$ with
\begin{gather*}
q=e^\hbar,\qquad q_i=q^{(\root{i},\root{i})/2}\\
[n]_i=\frac{q_i^n-q_i^{-n}}{q_i-q_i^{-1}}\aand
[n]_i!=[n]_i[n-1]_i\cdots[1]_i
\end{gather*}

\subsection{Monodromy conjecture}

It is known that the Yangian $Y(\g)$ and the quantum loop algebra
$\Uhghat$ have the same \fd representation theory (see \cite{Va}
and \cite{GTL1}). By analogy with the \qW description of the monodromy
of the (rational) Casimir connection of $\g$ conjectured by De Concini
(unpublished) and independently in \cite{TL1,TL2}, and proved in
\cite{TL1,TL3}, we make the following

\begin{conjecture}
The monodromy of the trigonometric Casimir connection is
equivalent to the quantum Weyl group action of the affine
braid group $\wh{B}$ on \fd $U_\hbar(L\g)$--modules.
\end{conjecture}

We will return to this conjecture in forthcoming work in collaboration
with S. Gautam \cite{GTL2}.

\section{The trigonometric Casimir connection of $\gl{n}$}\label{se:trigo gl}

We consider in this section the Yangian $Y(\gl{n})$ of the Lie algebra $\gl{n}$.
The latter does not possess a presentation of the form given in \S \ref{ss:J(x)}
but may be defined via a ternary, or $RTT$ presentation. By exploiting the
interplay between the latter and its loop presentation, we construct a flat,
trigonometric connection with values in $Y(\gl{n})$. We then relate it to the
corresponding connection for $\sl{n}$ and show that, when it is taken with
values in a tensor product of evaluation modules, it coincides with the
trigonometric dynamical equations \cite{TV}.

\subsection{The $RTT$ presentation of $Y(\gl{n})$}

The Yangian $Y(\gl{n})$ is the unital, associative algebra over $\Ch$ generated
by elements $t_{ij}^{(r)}$, $1\leq i,j\leq n$, $r\geq 1$, subject to the relations
\footnote{we follow here the conventions of \cite{Mo} and \cite{NO}.}
\begin{equation}\label{eq:tttt}
[t_{ij}^{(r+1)},t_{kl}^{(s)}]-[t_{ij}^{(r)},t_{kl}^{(s+1)}]=
t_{kj}^{(r)}t_{il}^{(s)}-t_{kj}^{(s)}t_{il}^{(r)}
\end{equation}
where $r,s\in\IN$ and $t_{ij}^{(0)}=\delta_{ij}$.

Let $V=\IC^n$ with standard basis $e_1,\ldots,e_n$ and let $E_{ij}e_k=\delta_{jk}
e_i$ be the corresponding basis of elementary matrices of $\gl{n}$. The map $
\imath:E_{ij}\to t_{ij}^{(1)}$ defines an embedding of $\gl{n}$ into $Y(\gl{n})$ and
we will often identify $\gl{n}$ with its image under $\imath$. Moreover, for every
$s\geq 1$, the subspace spanned by the elements $t_{ij}^{(s)}$ transforms like
the adjoint representation under the commutator action of $\gl{n}$.

The above relations may be more compactly rewritten as follows. For any
$r\geq 0$, let $T^{(r)}$ be the $n\times n$ matrix with values in $Y(\gl{n})$
given by
$$T^{(r)}=\sum_{1\leq i,j\leq n} E^V_{ij}\otimes t_{ij}^{(r)}$$
where $E_{ij}^V$ are again elementary matrices and the superscript $V$
is used to stress the fact that they should be thought of as elements of the
algebra $\End(V)$ rather than the underlying Lie algebra $\gl{n}$.

Let $u$ be a formal variable and set
$$T=\sum_{r\geq 0}T^{(r)}u^{-r}\in \End(V)\otimes Y(\gl{n})\fmlu$$
Finally, let
$$R(u)=1-Pu^{-1}\in \End(V\otimes V)\fmlu$$
be Yang's $R$--matrix, where $P\in\End(V\otimes V)$ acts as the permutation
of the two tensor factors. Then the relations \eqref{eq:tttt} are equivalent to 
$$R(u-v)T_1(u)T_2(v)=T_2(v)T_1(u)R(u-v)$$
where
$T_1(u),T_2(u)\in\End(V\otimes V)\otimes Y(\gl{n})\fmlu$
are given by
$$T_1(u)=\sum_{i,j,r}E^V_{ij}\otimes 1\otimes t_{ij}^{(r)}\,u^{-r}\aand
T_2(v)=\sum_{i,j,r}1\otimes E^V_{ij}\otimes t_{ij}^{(r)}\,v^{-r}$$

\subsection{The loop presentation of $Y(\gl{n})$}

Let $E(u),H(u),F(u)$ be the factors of the Gauss decomposition of $T(u)$.
Specifically,
\begin{gather*}
F(u)=1+\sum_{i>j}E^V_{ij}\otimes f_{ij}(u)\qquad
E(u)=1+\sum_{i<j}E^V_{ij}\otimes e_{ij}(u)\\
H(u)=\sum_i E^V_{ii}\otimes h_i(u),\quad
\end{gather*}
are, respectively, the unique lower unipotent, upper unipotent and diagonal
matrices with coefficients in $Y(\gl{n})\fmlu$ such that
\begin{equation}\label{eq:Gauss}
T(u)=F(u)H(u)E(u)
\end{equation}
Noting that $H(u),E(u),F(u)=1\mod u^{-1}$, write
$$h_i(u)=1+\sum_{r\geq 1}h_i^{(r)}u^{-r},\quad
f_{ij}(u)=\sum_{r\geq 1}f_{ij}^{(r)}u^{-r},\quad
e_{ij}(u)=\sum_{r\geq 1}e_{ij}^{(r)}u^{-r}$$
The coefficients of $e_{ii+1}(u),f_{ii+1}(u)$ and $h_i(u)$ give another system of
generators of $Y(\gl{n})$. Moreover, The elements $h_i(u)$ commute and their
coefficients generate a maximal commutative subalgebra of $Y(\gl{n})$ called
the {\it Gelfand--Zetlin subalgebra} $H_n$.

\subsection{}\label{ss:h t}

The Gauss decomposition \eqref{eq:Gauss} yields in particular
$$t_{ij}^{(1)}=\left\{\begin{array}{cl}
e_{ij}^{(1)}&\text{if $i<j$}\\[.2em]
h_i^{(1)}&\text{if $i=j$}\\[.2em]
f_{ij}^{(1)}&\text{if $i>j$}
\end{array}\right.$$
which we will use to identify the copies of $\gl{n}$ inside each presentation.
Moreover,
\begin{equation}\label{eq:h to t}
\begin{split}
t_{ii}^{(2)}
&=
h_i^{(2)}+\sum_{j<i}E_{ij}E_{ji}\\
&=
h_i^{(2)}+\half{1}\sum_{j<i}\left(\kappa_{\theta_j-\theta_i}-(E_{jj}-E_{ii})\right)
\end{split}
\end{equation}
where $\theta_a$ is the linear form given by $\theta_a(E_{bb})=\delta_{ab}$
and $\kappa_{\theta_a-\theta_b}=E_{ab}E_{ba}+E_{ba}E_{ab}$ is the truncated
Casimir operator corresponding to the root $\theta_a-\theta_b$.

\subsection{}

Define the elements $D_i\in H_n$ by
\begin{align}
D_i
&=
2h_i^{(2)}-\sum_{j<i}(E_{jj}-E_{ii})-E_{ii}^2
\label{eq:Di hi}\\
&=
2t_{ii}^{(2)}-\sum_{j<i}\kappa_{\theta_j-\theta_i}-E_{ii}^2
\label{eq:Di ti}
\end{align}
The symmetric group $\SS_n$ acts by algebra automorphisms on $Y(\gl{n})$ by
\begin{equation}\label{eq:Sn on Y}
\sigma(t_{ij}^{(r)})=t_{\sigma(i)\sigma(j)}^{(r)}
\end{equation}

\begin{lemma}\label{le:Di}
The following holds
\begin{enumerate}
\item $[D_i,D_j]=0$.
\item $(i\, i+1)D_j=D_j$ if $j\notin\{i,i+1\}$.
\item $(i\, i+1)D_i-D_{i+1}=\kappa_{\theta_i-\theta_{i+1}}$.
\end{enumerate}
\end{lemma}
\proof (1) follows from \eqref{eq:Di hi} and the fact that the $h_i^{(r)}$ commute.
(2) and (3) follows from \eqref{eq:Di ti}. \halmos

\subsection{}

The following is a direct consequence of \eqref{eq:Di ti} and \eqref{eq:h to t}

\begin{lemma}\label{le:D}
The element $\bfD=D_1+\cdots+D_n$ is given by
\begin{equation*}
\begin{split}
\bfD
&=2\sum_i t_{ii}^{(2)}-C_{\gl{n}}\\
&=2\sum_i h_i^{(2)}-2\rho^\vee-\sum_i E_{ii}^2
\end{split}
\end{equation*}
where
$$C_{\gl{n}}=\sum_{i<j}\kappa_{\theta_i-\theta_j}+\sum_i E_{ii}^2\aand
2\rho^\vee=\sum_{i<j}(E_{ii}-E_{jj})$$
are the Casimir operator and sum of the positive coroots of $\gl{n}$.
\end{lemma}

\subsection{The trigonometric Casimir connection of $\gl{n}$}\label{ss:trigo gl}

Let $D_i$ be given by \eqref{eq:Di ti} and define elements $\Delta_i\in H_n$
by (cf. \S \ref{ss:delta})
\begin{equation*}
\begin{split}
\Delta_i
&=
D_i-\half{1}\sum_{a<b}(\theta_a-\theta_b)(E_{ii})\,\kappa_{\theta_a-\theta_b}\\
&=
2t_{ii}^{(2)}-\half{1}\sum_{j\neq i}\kappa_{\theta_i-\theta_j}-E_{ii}^2
\end{split}
\end{equation*}
Let $H\subset GL_n$ be the maximal torus consisting of diagonal matrices,
$H\reg$ its set of regular elements and $\Y(\gl{n})$ the trivial $Y(\gl{n})
$--bundle over $H\reg$.

\begin{definition}
The {\it trigonometric Casimir connection} of $\gl{n}$ is the connection on
$\Y(\gl{n})$ given by either of the following forms
\begin{align*}
\nabla
&=d
-\sum_{i<j}\frac{d(\theta_i-\theta_j)}{e^{\theta_i-\theta_j}-1}\,
\kappa_{\theta_i-\theta_j}
-\sum_{i=1}^n d\theta_i\,D_i\\[1 em]
&=d
-\half{1}\sum_{i<j}\frac{e^{\theta_i-\theta_j}+1}{e^{\theta_i-\theta_j}-1}{d(\theta_i-\theta_j)}\,
\kappa_{\theta_i-\theta_j}
-\sum_{i=1}^n d\theta_i\,\Delta_i
\end{align*}
\end{definition}

\subsection{}

Let the symmetric group $\SS_n$ act on the vector bundle $\Y(\gl{n})$
by permutations of the base and automorphisms \eqref {eq:Sn on Y} of
the fibre.

\begin{theorem}\label{th:trigo gl}
The trigonometric Casimir connection of $\gl{n}$ is flat and equivariant
under $\SS_n$.
\end{theorem}
\proof Let $\ol{H},U\subset H$ be the subtori consisting respectively of
diagonal matrices of determinant 1 and multiples of the identity, and let
$\ol{\h},\u=\IC\bf1_n$ and $\h$ be their Lie algebras, where ${\bf1_n}=
\sum_{i=1}^nE_{ii}$. Thus
$$\h=\ol{\h}\oplus\IC\bf1_n\aand
\h^*\cong\ol{\h}^*\oplus\IC\tr$$
where $\tr:\h\to\IC$ is the trace. 
Clearly, $H\cong\ol{H}\times U$ and the connection $\nabla$ decomposes
as the product of the following $Y(\gl{n})$--valued connections on $\ol{H}
\reg$ and $U$ respectively
\begin{align}
\ol{\nabla}
&=
d-\sum_{i<j}\frac{d(\theta_i-\theta_j)}{e^{\theta_i-\theta_j}-1}\,
\kappa_{\theta_i-\theta_j}
-du_a\,D(u^a)
\label{eq:trigo gl restricted}\\
\nabla^U
&=d-\frac{1}{n}d\tr\,\bfD
\end{align}
where $D:\h\to Y(\gl{n})$ is given by $D(E_{ii})=D_i$ and $\{u_a\},\{u^a\}$
are dual bases of $\ol{\h}^*$ and $\ol{\h}$ respectively. $\nabla^U$ is clearly
flat and equivariant under the action of $\SS_n$ since the latter acts trivially
on $U$ and, by Lemma \ref{le:D} on $\bfD$. Since the latter commutes with
the coefficients of $\ol{\nabla}$, the flatness and equivariance of $\nabla$
reduces to that of $\ol{\nabla}$ which, in turn is determined by Theorem
\ref{th:trigo flat} and Proposition \ref{pr:equivariance}.
The relations \eqref{eq:quadratic} have already been checked in the proof
of the flatness of the trigonometric Casimir connection for $\sl{n}$ in Theorem
\ref{th:trigo Casimir}. The relations \eqref{eq:comm} and the equivariance
relations \eqref{eq:equiv2} follow from Lemma \ref{le:Di}. There remains to check
that, for any $i=1,\ldots,n-1$ and $u\in\Ker(\theta_i-\theta_{i+1})$, $[\kappa_
{\theta_i-\theta_{i+1}},D(u)]=0$. This reduces to checking that $[\kappa_{\theta_i
-\theta_{i+1}},D_j]=0$ for $j\notin\{i,i+1\}$ and that $[\kappa_{\theta_i
-\theta_{i+1}},D_i+D_{i+1}]=0$ which follows easily from \eqref{eq:Di ti}. \halmos
 
\subsection{Rational form}

It is well known that trigonometric connections of type $GL_n$ may be put into
a rational form, thus expressing them as \KZ type connections on $n+1$ points,
with one frozen to $0$. We carry this step below for the connection \ref{ss:trigo gl}.

Let $z_i=e^{\theta_i}$, $i=1,\ldots,n$ be the standard coordinates on the
torus $H\cong(\IC^*)^n$ of $GL_n$.\footnote{note that these differ from
the coordinates $z_i=e^{\alpha_i}=e^{\theta_i-\theta_{i+1}}$ used in \S
\ref {ss:compactification}.} Since $d\theta_i=dz_i/z_i$ and
$$\frac{d(\theta_i-\theta_j)}{e^{\theta_i-\theta_j}-1}=
\frac{d(z_i-z_j)}{z_i-z_j}-\frac{dz_i}{z_i}$$
the connection \ref{ss:trigo gl} is equal to
$$\nabla=d-\sum_{i<j}\frac{d(z_i-z_j)}{z_i-z_j}\,
\kappa_{\theta_i-\theta_j}-
\sum_i\frac{dz_i}{z_i}\,\wt{D}_i$$
where
$$\wt{D_i}=
D_i-\sum_{j>i}\kappa_{\theta_i-\theta_j}=
2t_{ii}^{(2)}-\sum_{j\neq i}\kappa_{\theta_i-\theta_j}-E_{ii}^2$$

\subsection{Relation to $Y(\sl{n})$}

Following Olshanski and Drinfeld, we realise the Yangian $Y(\sl{n})$ as a Hopf
subalgebra of $Y(\gl{n})$ as follows (see \cite[\S 1.8]{Mo}). Let $\IA=1+u^{-1}\IC
\fmlu$ be the abelian group of formal power series in $u^{-1}$ with constant term
1. $\IA\ni f$ acts on $Y(\gl{n})$ by Hopf algebra automorphisms given by
$$T(u)\to f(u)T(u)$$
The Hopf subalgebra $Y(\gl{n})^{\IA}\subset Y(\gl{n})$ of elements fixed by $\IA$
is isomorphic to $Y(\sl{n})$.

\subsection{}

The generators $T_{i,r}$ of the presentation of $Y(\sl{n})$ described
in \S \ref{ss:new realisation} may be obtained within the RTT realisation of
$Y(\gl{n})$ as follows \cite[Rk. 5.12]{BK}. Consider their generating
function $T_i(u)=1+\sum_{r\geq 0}T_{i,r}u^{-r-1}$. Then,
$$T_i(u)=h_i(u-\frac{i-1}{2})^{-1}\cdot h_{i+1}(u-\frac{i-1}{2})$$
To spell this out, consider a formal power series $a(u)=1+a_1u^{-1}+a_2u
^{-2}+\cdots$. Then, for any $\lambda \in\IC$ one has
\begin{equation*}
\begin{split}
a(u-\lambda)
&=
1+a_1u^{-1}(1-\frac{\lambda}{u})^{-1}+a_2u^{-2}(1-\frac{\lambda}{u})^{-2}+
\cdots\\
&=
1+a_1u^{-1}+(a_2+\lambda a_1)u^{-2}+\cdots
\end{split}
\end{equation*}
and therefore
$$a(u-\lambda)^{-1}=1-a_1u^{-1}-(a_2+\lambda a_1-a_1^2)u^{-2}+\cdots$$
It follows that $T_{i,0}=-(E_{ii}-E_{i+1,i+1})$ and
\begin{equation}\label{eq:kappa h}
T_{i,1}=
-(h_i^{(2)}-h_{i+1}^{(2)})-\frac{i-1}{2}(E_{ii}-E_{i+1,i+1})+E_{ii}^2-E_{ii}E_{i+1,i+1}
\end{equation}

\subsection{}

Let $\ol{H}\subset H$ be the torus of $SL_n$ consisting of diagonal matrices
with determinant 1.

\begin{proposition}
The restriction of the trigonometric Casimir connection of $\gl{n}$ to $\ol
{H}\reg$ takes values in $Y(\sl{n})$ and is equal to the sum of the trigonometric
Casimir connection of $\sl{n}$ with the $\ol{\h}\subset Y(\sl{n})$--valued,
closed one--form
$$-\sum_{i=1}^n d\lambda_i\,(E_{ii}-E_{i+1,i+1})$$
where $\{\lambda_i\}$ are the fundamental weights of $\sl{n}$.
\end{proposition}
\proof The restriction of the trigonometric Casimir connecton of $\gl{n}$ to
$\ol{H}\reg$ is given by \eqref{eq:trigo gl restricted}, namely
$$\ol{\nabla}=
d-\sum_{i<j}\frac{d(\theta_i-\theta_j)}{e^{\theta_i-\theta_j}-1}\,
\kappa_{\theta_i-\theta_j}
-du_i\,D(u^i)$$
where $D:\h\to Y(\gl{n})$ is given by $D(E_{ii})=D_i$ and $\{u_i\},\{u^i\}$ are
dual bases of $\ol{\h}^*$ and $\ol{\h}$ respectively. Choosing $u_i=\lambda_i$
so that $u^i=E_{ii}-E_{i+1,i+1}$, $i=1,\ldots,n-1$ and comparing with the form
\eqref{eq:T tau form} we need to show that
$$D_i-D_{i+1}=-2T_{i,1}+(E_{ii}-E_{i+1,i+1})^2+(E_{ii}-E_{i+1,i+1})$$
By \eqref{eq:Di hi}, the \lhs is equal to
$$2(h_i^{(2)}-h_{i+1}^{(2)})+i\,(E_{ii}-E_{i+1,i+1})-E_{ii}^2+E_{i+1,i+1}^2$$
and the result follows from \eqref{eq:kappa h} \halmos

\subsection{Evaluation homomorphism}

The Yangian $Y(\gl{n})$ possesses an evaluation homomorphism $\ev:
Y(\gl{n})\to U\gl{n}$ defined by
$$\ev(t_{ij}(u))=\delta_{ij}+E_{ij}u^{-1}$$
where $t_{ij}(u)=\sum_{r\geq 0}t_{ij}^{(r)}u^{-r}$.
When composed with the translation automorphisms $\tau_a$, $a\in\IC$
given by $\tau_a T(u)=T(u-a)$, that is
$$\tau_a T^{(r)}=
\delta_{r0}+
\sum_{s=1}^r T^{(s)}\begin{pmatrix}r-1\\r-s\end{pmatrix}a^{r-s}$$
this yields a one--parameter family of evaluation homomorphisms $\ev_a
=\ev\circ\tau_a$ given by
\begin{equation}\label{eq:ev tij}
\ev_a(t_{ij}^{(r)})=\delta_{r0}\delta_{ij}+E_{ij}a^{r-1}
\end{equation}

\subsection{Hopf algebra structure}

$Y(\gl{n})$ is a Hopf algebra with coproduct
$$\Delta(t_{ij}(u))=\sum_{k=1}^n t_{ik}(u)\otimes t_{kj}(u)$$
For any $m\geq 2$, let $\Delta^{(m)}:Y(\gl{n})\to Y(\gl{n})^{\otimes m}$
denote the corresponding iterated coproduct. Then,
\begin{equation}\label{eq:Delta tij}
\Delta^{(m)}(t_{ij}^{(2)})=
\sum_{p=1}^m (t_{ij}^{(2)})_p+
\sum_{\substack{1\leq k\leq n\\[.3ex]1\leq p<q\leq m}}(E_{ik})_p(E_{kj})_q
\end{equation}
where $X_p=1^{\otimes (p-1)}\otimes X\otimes 1^{\otimes (m-p)}$.

\subsection{Evaluation modules}

For any $\ul{a}=(a_1,\ldots,a_m)\in\IC^m$ define $\ev_{\ul{a}}:Y(\gl{n})\to
U\gl{n}^{\otimes m}$ by
$$\ev_{\ul{a}}=\ev_{a_1}\otimes\cdots\otimes\ev_{a_m}\circ\Delta^{(m)}$$

\begin{proposition}\label{pr:ev trigo Casimir gl}
The image of the trigonometric Casimir connection of $\gl{n}$ under the
homomorphism $\ev_{\ul{a}}$ is the $U\gl{n}^{\otimes m}$--valued connection
given by
$$\nabla_{\ul{a}}=d
-\sum_{i<j}\frac{d(\theta_i-\theta_j)}{e^{\theta_i-\theta_j}-1}\,
\Delta^{(m)}(\kappa_{\theta_i-\theta_j})
-\sum_{i=1}^n d\theta_i\,D_{i,\ul{a}}$$
where
\begin{equation*}
D_{i,\ul{a}}=
2\sum_{p=1}^m a_p (E_{ii})_p
+2\negthickspace\negthickspace\negthickspace\negthickspace
\sum_{\substack{1\leq j\leq n\\[.3ex]1\leq p<q\leq m}}(E_{ij})_p(E_{ji})_q
-\sum_{j<i}\Delta^{(m)}(\kappa_{\theta_j-\theta_i})
-\Delta^{(m)}(E_{ii}^2)
\end{equation*}
\end{proposition}
\proof By construction $D_{i,\ul{a}}=\ev_{\ul{a}}(D_i)$ and is given by the
above expression by \eqref{eq:Di ti}, \eqref{eq:Delta tij} and \eqref{eq:ev tij}.
\halmos

\subsection{The trigonometric dynamical differential equations for $\gl{n}$}

In \cite{TV}, Tarasov and Varchenko considered differential operators $\D_1,
\ldots,\D_n$ in the variables $z_1,\ldots,z_n\in\IC^\times$ with coefficients in
$U\gl{n}^{\otimes m}$ given by
$$\D_i=z_i\partial_{z_i}+\lambda L_i(\ul{a},\ul{z})$$
where $\lambda\in\IC$,
$\ul{a}=(a_1,\ldots,a_m)\in\IC^m$, $\ul{z}=(z_1,\ldots,z_n)$ and
\begin{multline*}
L_i(\ul{a},\ul{z})
=\half{\Delta^{(m)}(E_{ii}^2)}-\sum_{p=1}^m a_p (E_{ii})_p\\
-\sum_{\substack{1\leq j\leq n\\[.3ex]1\leq p<q\leq m}}(E_{ij})_p(E_{jk})_q
-\sum_{j\neq i}\frac{z_j}{z_i-z_j}\Delta^{(m)}(E_{ij}E_{ji}-E_{ii})
\end{multline*}
Set $z_i=e^{\theta_i}$ so that $z_i\partial_{z_i}=\partial_{\theta_i}$ and
$$\frac{z_j}{z_i-z_j}=
\frac{1}{e^{\theta_i-\theta_j}-1}=-\left(\frac{1}{e^{\theta_j-\theta_i}-1}+1\right)$$
Since
$$E_{ij}E_{ji}-E_{ii}=\half{1}(\kappa_{\theta_i-\theta_j}-(E_{ii}+E_{jj}))$$
the operators $\D_i$ are the covariant derivatives for the connection
$$\nabla'_{\ul{a}}=d-\half{\lambda}\left(\sum_{i<j}\frac{d(\theta_i-\theta_j)}{e^{\theta_i-\theta_j}-1}
\Delta^{(m)}(\kappa_{\theta_i-\theta_j}-(E_{ii}+E_{jj}))
+\sum_{i=1}^n d\theta_i\, D_{i,\ul{a}}'\right)$$
where
\begin{multline*}
D_{i,\ul{a}}'=
-\Delta^{(m)}(E_{ii}^2)
+2\sum_{p=1}^m a_p (E_{ii})_p
+2\negthickspace\negthickspace\negthickspace\negthickspace
\sum_{\substack{1\leq j\leq n\\[.3ex]1\leq p<q\leq m}}(E_{ij})_p(E_{ji})_q\\
-\sum_{j<i}\Delta^{(m)}(\kappa_{\theta_j-\theta_i}-(E_{ii}+E_{jj}))$$
\end{multline*}
By Proposition \ref{pr:ev trigo Casimir gl}, $\nabla'_{\ul{a}}$ is the image of the
trigonometric Casimir connection for $\gl{n}$\footnote{when the latter is scaled
by a factor of $\lambda/2$, as in Remark \ref{rk:scaling}.} under the homomorphism
$\ev_{\ul{a}}:Y(\gl{n})\to (U\gl{n})^{\otimes m}$ plus the $\h$--valued, closed one--form
$$\half{\lambda}\,\Delta^{(m)}\left(\sum_{i<j}\frac{d(\theta_i-\theta_j)}{e^{\theta_i-\theta_j}-1}
(E_{ii}+E_{jj})
-\sum_id\theta_i
\sum_{j<i}
(E_{ii}+E_{jj})\right)$$

\section{Bispectrality}\label{se:bispesctral}

We show in this section that the trigonometric Casimir connection with values
in a tensor product of $Y(\g)$--modules commutes with the $q$KZ difference
equations of \FR determined by the rational $R$--matrix of $Y(\g)$. This was
checked by \TV for $\g=\gl{n}$ in the case where all representations are evaluation
modules \cite{TV}.

\subsection{Hopf algebra structure \cite{Dr1}}

If $\g$ is simple, $Y(\g)$ is a Hopf algebra with coproduct $\Delta:
Y(\g)\to Y(\g)\otimes Y(\g)$ given on generators by
\begin{align*}
\Delta(x)&=x\otimes 1+1\otimes x\\
\Delta(J(x))&=J(x)\otimes 1+1\otimes J(x)+\half{\hbar}[x\otimes 1,t]
\end{align*}
where $t=\sum_a x_a\otimes x^a\in(\g\otimes\g)^\g$, with $\{x_a\},
\{x^a\}$ dual bases of $\g$ \wrt the given inner product. Thus, if $
\Delta^{(n)}:Y(\g)\to Y(\g)^{\otimes n}$ is the iterated coproduct,
then
\begin{align}
\Delta^{(n)}(x)&=\sum_{i=1}^n x^{(i)}\\
\Delta^{(n)}(J(x))&=
\sum_{i=1}^n J(x)^{(i)}+\frac{\hbar}{2}\sum_{1\leq i<j\leq n}[x^{(i)},t^{ij}]
\label{eq:Delta k}
\end{align}
where $x^{(i)}=1^{\otimes(i-1)}\otimes x\otimes 1^{\otimes (n-i)}$ and
$t^{ij}=\sum_a x_a^{(i)}{(x^a)}^{(j)}$.

\subsection{Translation automorphisms \cite{Dr1}}

$Y(\g)$ possesses a one--parameter group of Hopf algebra automorphisms
$T_v$, $v\in\IC$ given by
$$T_v x=x\aand T_v J(x)=J(x)+vx$$
If $v_1,\ldots,v_n\in\IC$, we set $T_{v_1,\ldots,v_n}=T_{v_1}\otimes\cdots
\otimes T_{v_n}\in\Aut(Y(\g)^{\otimes n})$ and
\begin{equation}\label{eq:translated delta}
\Delta_{v_1,\ldots,v_n}=T_{v_1,\ldots,v_n}\circ\Delta^{(n)}:Y(\g)\to Y(\g)^{\otimes n}
\end{equation}

\subsection{The universal $R$--matrix of $Y(\g)$ \cite{Dr1}}

Let $R(u)\in Y(\g)\otimes Y(\g)[[u^{-1}]]$ be the universal $R$--matrix of $Y(\g)$.
$R(u)$ satisfies
\begin{gather}
\Delta\otimes\id (R(u))=R^{13}(u)R^{23}(u)
\label{eq:cabling 1}\\
\id\otimes\Delta (R(u))=R^{13}(u)R^{12}(u)
\label{eq:cabling 2}\\
\Delta^{21}(x)=R(u)\Delta(x)R(u)^{-1}
\label{eq:almost coco}\\
T_{v,w}R(u)=R(u+v-w)
\label{eq:translation}
\end{gather}
The above relations imply that $R$ satisfies the \QYBE (QYBE) with spectral
parameter
$$R^{12}(u)R^{13}(u+v)R^{23}(v)=R^{23}(v)R^{13}(u+v)R^{12}(u)$$
and the more general form of \eqref{eq:almost coco}
$$R(u)(T_{v,w}\Delta(x))R(u)^{-1}=T_{v,w}\Delta^{21}(x)$$

\subsection{The rational qKZ equations \cite{FR}}

Let $V_1,\ldots,V_n$ be $Y(\g)$--modules and $d_i\in GL(V_i)$ be such that
$d_id_jR^{ij}(u)=R^{ij}(u)d_id_j$ for any $1\leq i<j\leq n$. Fix a step $\varkappa
\in\IC^\times$, let $a_1,\ldots,a_n\in\IC$ be distinct and define operators
$$A_i=A_i(a_1,\ldots,a_n)\in\End(V_1\otimes\cdots\otimes V_n)$$
for $i=1,\ldots,n$ by
\begin{equation*}
\begin{split}
A_i
&=
R^{i-1\,i}(a_{i-1}-a_i-\varkappa)^{-1}R^{i-2\,i}(a_{i-2}-a_i-\varkappa)^{-1}\cdots
R^{1\,i}(a_1-a_i-\varkappa)^{-1}\cdot \\
&\phantom{==}
d_i\cdot
R^{i\,n}(a_i-a_n)R^{i\,n-1}(a_i-a_{n-1})\cdots R^{i\,i+1}(a_i-a_{i+1})
\end{split}
\end{equation*}
The (rational) $q$KZ equations of \FR are the system of difference equations
$T_i f=A_i f$ where $f$ takes values in $V_1\otimes\cdots\otimes V_n$ and
$$T_i f(a_1,\ldots,a_n)=f(a_1,\ldots,a_{i-1},a_i+\varkappa,a_{i+1},\ldots,a_n)$$
They are a consistent system in that $[\TT_i,\TT_j]=0$, where $\TT_i=A_i^{-1}T_i$
are the covariant difference operators.

\subsection{}

\begin{lemma}
The following holds for any $i=1,\ldots,n$
$$\TT_1\TT_2\cdots\TT_i=
(\wt{A}_i)^{-1}T_1\cdots T_i$$
where $\wt{A}_i=\wt{A}_i(a_1,\ldots,a_n)$ is given by
\begin{multline}\label{eq:A_i 1}
\wt{A}_i
=
d_1\cdots d_i
\left(R^{1\,n}(a_1-a_n)\cdots R^{i\,n}(a_i-a_n)\right)\\
\left(R^{1\,n-1}(a_1-a_{n-1})\cdots R^{i\,n-1}(a_i-a_{n-1})\right)\cdots\\
\cdots
(R^{1\,i+1}(a_1-a_{i+1})\cdots R^{i\,i+1}(a_i-a_{i+1}))
\end{multline}
Thus, $\wt{A}_n=d_1\cdots d_n$ and for any $i\leq n-1$
\begin{equation}\label{eq:A_i 2}
\wt{A}_i
=
d_1\cdots d_i\,
\Delta^{(i)}_{a_1-a_i,\ldots,a_{i-1}-a_i,0}\otimes\Delta^{(n-i)}_{a_{i+1}-a_n,\ldots,a_{n-1}-a_n,0}
\left(R(a_i-a_n)\right)
\end{equation}
\end{lemma}
\proof Clearly,
$$\TT_1\TT_2\cdots\TT_i=
\left(T_1\cdots T_{i-1}(A_i)\cdots T_1(A_2)A_1\right)^{-1}T_1\cdots T_i$$
and, for any $j\leq i$,
$$T_1\cdots T_{j-1}(A_j)=
(R^{i-1\,i})^{-1}(R^{i-2\,i})^{-1}\cdots(R^{1\,i})^{-1}
\cdot d_i\cdot
R^{i\,n}R^{i\,n-1}\cdots R^{i\,i+1}
$$
where $R^{kl}$ is shorthand for $R^{kl}(a_k-a_l)$. The first claimed identity now
follows by induction on $i$ using the QYBE. The second follows from the relations
\eqref{eq:cabling 1}--\eqref{eq:cabling 2} and \eqref{eq:translation}. \halmos

\subsection{Bispectrality}

To couple the $q$KZ and trigonometric Casimir connection equations with values
in the tensor product $V_1\otimes\cdots\otimes V_n$, assume that each $V_i$
is an integrable $\g$--module and that $d_i$ is the $GL(V_i)$--valued function
on the torus $H$ given by
\begin{equation}\label{eq:di}
d_i(e^u)=(e^{-u})^{(i)}
\end{equation}
Let also $\nabla_{\ul{a}}$ be the trigonometric Casimir connection with values
in the $Y(\g)$--module $T_{-a_1,\ldots,-a_n}^*\,V_1\otimes\cdots\otimes V_n$
and scaled by a factor of $2\varkappa$ as in Remark \ref{rk:scaling}. Thus, $\nabla
_{\ul{a}}$ is the $\End(V_1\otimes\cdots\otimes V_n)$--valued connection given
by
$$\nabla_{\ul{a}}=d-\frac{1}{2\varkappa}\,\Delta^{(n)}_{a_1,\ldots,a_n}(B)$$
where
$$B=\frac{\hbar}{2}
\sum_{\alpha\in\Phi_+}\frac{e^\alpha+1}{e^\alpha-1}d\alpha\, \kappa_\alpha
-2du_i\,J(u^i)$$

\begin{theorem}\label{th:bispectral}
The $qKZ$ operators $\TT_i$ commute with the trigonometric
Casimir connection $\nabla_{\ul{a}}$.
\end{theorem}
\proof It suffices to prove that $[\nabla_{\ul{a}},\TT_1\cdots\TT_i]=0$ for any $i=1,
\ldots,n$.
Since
\begin{multline*}
[d-(2\varkappa)^{-1}\Delta^{(n)}_{a_1,\ldots,a_n}(B),
(\wt{A}_i)^{-1}T_1\cdots T_i]\,(T_1\cdots T_i)^{-1}\\[1.1ex]
=
d\wt{A}_i^{-1}
-
(2\varkappa)^{-1}\wt{A}_i^{-1}(\id-\Ad(T_1\cdots T_i))\Delta^{(n)}_{a_1,\ldots,a_n}(B)\\
-
(2\varkappa)^{-1}[\Delta^{(n)}_{a_1,\ldots,a_n}(B),(\wt{A}_i)^{-1}]
\end{multline*}
the claim follows from the two lemmas below. \halmos

\subsection{}

\begin{lemma}
\begin{equation}\label{eq:cross diff}
(d\wt{A}_i)\wt{A}_i^{-1}=
(2\varkappa)^{-1}(\Ad(T_1\cdots T_i)-\id)\Delta^{(n)}_{a_1,\ldots,a_n}(B)
\end{equation}
\end{lemma}
\proof By \eqref{eq:A_i 1}, the \lhs of \eqref{eq:cross diff} is equal to
$$d(d_1\cdots d_i)\,(d_1\cdots d_i)^{-1}=-\sum_{j=1}^i du_a (u^a)^{(j)}$$
where we used \eqref{eq:di}. Write $B=B_1+B_2$ where
$$B_1=\frac{\hbar}{2}
\sum_{\alpha\in\Phi_+}d\alpha\, \frac{e^\alpha+1}{e^\alpha-1}\kappa_\alpha
\aand
B_2=-2du_i\,J(u^i)$$
Since $B_1$ takes values in $U\g$, $\Delta_{a_1,\ldots,a_n}(B_1)^{(n)}$ is 
independent of $a_1,\ldots,a_n$ and the \rhs of \eqref{eq:cross diff} is equal
to $(2\varkappa)^{-1}(\Ad(T_1\cdots T_i)-\id)\Delta^{(n)}_{a_1,\ldots,a_n}(B_2)$.
By \eqref{eq:Delta k}, for any $x\in\g$,
$$(\Ad(T_1\cdots T_i)-\id)\Delta^{(n)}_{a_1,\ldots,a_n}(J(x))=
\sum_{j=1}^i \varkappa x^{(i)}$$
so that the \rhs of \eqref{eq:cross diff} is equal to $-\sum_{j=1}^i du_a(u^a)^{(j)}$.
\halmos

\subsection{}

\begin{lemma}
$$[\wt{A}_i,\Delta^{(n)}_{a_1,\ldots,a_n}(B)]=0$$
\end{lemma}
\proof For any $x\in Y(\g)$ and $1\leq i\leq n$,
$$\Delta_{a_1,\ldots,a_n}^{(n)}(x)=
\Delta_{a_1-a_i,\ldots,a_{i-1}-a_i,0}^{(i)}\otimes\Delta^{(n-i)}_{a_{i+1}-a_n,\ldots,a_{n-1}-a_n,0}
\circ\Delta_{a_i,a_n}(x)$$
so that, by \eqref{eq:A_i 2} and the fact that $d_1\cdots d_i=\Delta^{(i)}(d_1)$, it suffices
to prove the claimed identity for $n=2$ and $i=1$. We have
\begin{multline}\label{eq:we have}
d_1^{-1}\,[d_1R(a_1-a_2),\Delta_{a_1,a_2}(B)]\,R(a_1-a_2)^{-1}\\
=(\id-\Ad(d_1^{-1}))\Delta_{a_1,a_2}(B)+
(\Ad(R(a_1-a_2))-\id)\Delta_{a_1,a_2}(B)
\end{multline}
Let $u\in\h$. By \eqref{eq:Delta k},
$$(\id-\Ad(d_1^{-1}))\Delta_{a_1,a_2}(J(u))=\half{\hbar}(\id-\Ad(d_1^{-1}))[u^{(1)},t]$$
and, by \eqref{eq:almost coco},
$$(\Ad(R(a_1-a_2))-\id)\Delta_{a_1,a_2}(J(u))=\half{\hbar}\,[u^{(2)}-u^{(1)},t]=-\hbar\,[u^{(1)},t]$$
Thus, the \rhs of \eqref{eq:we have} with $B$ replaced by $B_2=-2du_a\,J(u^a)$ is
equal to
$$\hbar\,du_a(\id+\Ad(d_1^{-1}))[{u^a}^{(1)},t]$$
Since $R(a_1-a_2)$ commutes with $\Delta_{a_1,a_2}(B_1)$, we have left to compute
$$(\id-\Ad(d_1^{-1}))\Delta_{a_1,a_2}(B_1)=
\hbar\sum_{\alpha\in\Phi_+}d\alpha\,
\frac{e^\alpha+1}{e^\alpha-1}(\id-\Ad(d_1^{-1}))\ol{t}_\alpha$$
where $\ol{t}_\alpha=x_\alpha\otimes x_{-\alpha}+x_{-\alpha}\otimes x_\alpha$, with
$x_{\pm\alpha}\in\g_{\pm\alpha}$ such that $(x_\alpha,x_{-\alpha})=1$, so that 
$\Delta(\kappa_\alpha)=\kappa_\alpha\otimes 1+1\otimes \kappa_\alpha+2\ol{t}_\alpha$. By \eqref{eq:di},
$$(\id-\Ad(d_1^{-1}))\ol{t}_\alpha=
(1-e^\alpha)(x_\alpha\otimes x_{-\alpha}-e^{-\alpha}x_{-\alpha}\otimes x_\alpha)$$
so that, for any $\alpha\in\Phi_+$ and $u\in\h$,
\begin{multline*}
\alpha(u)
\frac{e^\alpha+1}{e^\alpha-1}(\id-\Ad(d_1^{-1}))\ol{t}_\alpha\\
=-\alpha(u)\left((e^\alpha+1)x_\alpha\otimes x_{-\alpha}-(e^{-\alpha}+1)x_{-\alpha}\otimes x_\alpha\right)\\
=-[u^{(1)},(\id+\Ad(d_1^{-1}))\ol{t}_\alpha]
\end{multline*}
whence the claimed result. \halmos

\begin{remark}
The proof of Theorem \ref{th:bispectral} works almost {\it verbatim} for
$\g=\gl{n}$ and gives the commutation of the rational $q$KZ connection
and trigonometric Casimir connections for $Y(\gl{n})$.
\end{remark}

\section{The affine KZ connection}\label{se:AKZ}

We show in this section that the \daha $\H'$ of $W$ is, very roughly speaking,
the 'Weyl group' of the Yangian $Y(\g)$. More precisely, we show that if $V$
is a $Y(\g)$--module whose restriction to $\g$ is {\it small}, the canonical
action of $W$ on the zero weight space $V[0]$ extends to one
of $\H'$. Moreover, the trigonometric Casimir connection with values in
$V[0]$ coincides with Cherednik's affine KZ connection with values
in this $\H'$--module.

\subsection{The degenerate affine Hecke algebra}

Let $K$ be the vector space of $W$--invariant functions $\Phi\to\IC$
and denote the natural linear coordinates on $K$ by $k_\alpha$, $\alpha
\in\Phi/W$. Recall \cite{Lu} that the degenerate affine Hecke algebra $\H'$
associated to $W$ is the algebra over $\IC[K]$ generated by the group
algebra $\IC W$ and the symmetric algebra $S\h$ subject to the relations
\begin{equation}\label{eq:DAHA}
s_i x_u-x_{s_i(u)}s_i=k_{\alpha_i}\alpha_i(u)
\end{equation}
for any simple reflection $s_i\in W$ and linear generator $x_u$, $u\in\h$,
of $S\h$.

\subsection{The affine KZ connection}

The AKZ connection is the trigonometric, $\H'$--valued connection given
by
\begin{equation}\label{eq:AKZ}
\nabla=
d-\sum_{\alpha\in\Phi_+}\frac{d\alpha}{e^\alpha-1}\,k_\alpha s_\alpha
-du_i \, x_{u^i}
\end{equation}
where $\{u_i\},\,\{u^i\}$ are dual bases of $\h^*,\,\h$ respectively. This
connection was defined by Cherednik in \cite{Ch1,Ch2} and proved to
be flat and $W$--equivariant. This may also be obtained as a consequence
of Theorem \ref{th:W trigo flat}. Indeed, the relations \eqref{eq:quadratic},
with $t_\alpha=k_\alpha s_\alpha$ are easily verified and, as pointed out
by Cherednik, the remaining relations are precisely those defining $\H'$.

\begin{remark} The $\delta$--form \eqref{eq:delta form} of the AKZ 
connection corresponds to Drinfeld's presentation of $\H'$ in terms
of $\IC W$ and non--commuting elements $y_u$ which transform
like the reflection representation of $W$ (see \cite{Dr3} and \cite
{RS}).\footnote{I owe this observation to Pavel Etingof.}
Indeed, it is given by
$$\nabla=d-\half{1}\sum_{\alpha\in\Phi_+}
\frac{e^\alpha-1}{e^\alpha-1}{d\alpha}\,k_\alpha s_\alpha
-du_i \, y_{u^i}$$
where the elements $y_u$ are defined by \eqref{eq:delta} as
\begin{equation}\label{eq:y x}
y_u=x_u-\half{1}\sum_{\alpha\in\Phi_+}\alpha(u)k_\alpha s_\alpha
\end{equation}
and therefore satisfy $s_i\, y_u s_i=y_{s_i(u)}$ by Proposition \ref
{pr:equivariance}.
\end{remark}

\subsection{$W$--action on zero weight spaces of $\g$--modules}

Let $G$ be the complex, simply--connected Lie group with Lie algebra $\g$,
$H$ the maximal torus with Lie algebra $\h$ and $N(H)\subset G$ its normaliser.
If $V$ is an integrable $\g$--module, the action of $N(H)$ on $V$ permutes
the weight spaces compatibly with the action of $N(H)$ on $H$. In particular,
it acts on the zero weight space $V[0]$ and this action factors through $W=
N(H)/H$.

\subsection{Small $\g$--modules}

Recall that a $\g$--module $V$ is {\it small} if $2\alpha$ is not a weight of $V$
for any root $\alpha$ \cite{Br,Re,Re2}. If $V$ is a small $\g$--module with a 
non--trivial zero weight space $V[0]$, the restriction to $V[0]$ of the square
$e_\alpha^2$ of a raising operator maps to the weight space $V[2\alpha]$
and is therefore zero. This implies the following result \cite[Prop. 9.1]{TL2}

\begin{lemma}\label{le:V[0]}
If $V$ is an integrable, small $\g$--module, the following holds on the zero
weight space $V[0]$
$$\kappa_\alpha=(\alpha,\alpha)(1-s_\alpha)$$
where the \rhs refers to the action of the reflection $s_\alpha\in W$ on $V[0]$.
\end{lemma}

\subsection{} 

Let $\H'_\hbar$ be the degenerate affine Hecke algebra of $W$ with
parameters
\begin{equation}\label{eq:k alpha}
k_\alpha=-\hbar(\alpha,\alpha)
\end{equation}

\begin{theorem}\label{th:trigo=AKZ}
Let $V$ be a \fd $Y(\g)$--module whose restriction to $\g$ is small.
\begin{enumerate}
\item The canonical $W$--action on the zero weight space $V[0]$ together
with either of the equivalent assignments
\begin{align}
x_u&\to-2T(u)_1+\half{1}\sum_{\alpha\in\Phi_+}k_\alpha\alpha(u)
\label{eq:H' action}\\
y_u&\to -2 J(u)
\label{eq:H' action y}
\end{align}
yield an action of $\H'_\hbar$ on $V[0]$.
\item The trigonometric Casimir connection of $\g$ with values in $\End
(V[0])$ is equal to the sum of the AKZ connection with values in the $\H'
_\hbar $--module $V[0]$ and the scalar valued one--form
\begin{equation}\label{eq:A}
\A=
\half{1}\sum_{\alpha\in\Phi}k_\alpha\frac{d\alpha}{e^\alpha-1}
\end{equation}
\end{enumerate} 
\end{theorem}
\proof The trigonometric Casimir connection with values in $\End(V[0])$
reads, by Lemma \ref{le:V[0]}
\begin{equation*}
\begin{split}
\nabla
&=
d-\hbar\sum_{\alpha\in\Phi_+}\frac{d\alpha}{e^\alpha-1}\,
(\alpha,\alpha)(1-s_\alpha)
+2du_i\,T(u^i)_1\\
&=
d-\sum_{\alpha\in\Phi_+}\frac{d\alpha}{e^\alpha-1}\,
k_\alpha s_\alpha
+2du_i\,T(u^i)_1+\sum_{\alpha\in\Phi_+}k_\alpha\frac{d\alpha}{e^\alpha-1}
\end{split}
\end{equation*}
where the weights $k_\alpha$ are given by \eqref{eq:k alpha}.
By \eqref{eq:trick}
$$\sum_{\alpha\in\Phi}(\alpha,\alpha)\frac{d\alpha}{e^\alpha-1}=
2\sum_{\alpha\in\Phi_+}(\alpha,\alpha)\frac{d\alpha}{e^\alpha-1}+
\sum_{\alpha\in\Phi_+}(\alpha,\alpha)d\alpha$$
so that if $\A$ is given by \eqref{eq:A}, then
$$\nabla-\A=
d-\sum_{\alpha\in\Phi_+}k_\alpha\frac{d\alpha}{e^\alpha-1}\,s_\alpha
+du_i\left(2T(u^i)_1-\half{1}\sum_{\alpha\in\Phi_+}k_\alpha\alpha(u^i)\right)$$
Applying Proposition \ref{pr:equivariance} to $\nabla-\A$
which is $W$--equivariant since $\nabla$ and $\A$ are, shows that the
map \eqref{eq:H' action} gives an action of $\H'_\hbar$ on $V[0]$. This
proves (1) and (2). The equivalence of \eqref{eq:H' action} and \eqref
{eq:H' action y} follows easily from \S \ref{ss:old new} and \eqref{eq:y x}.
\halmos

\begin{remark}
Theorem \ref{th:trigo=AKZ} extends to the trigonometric setting the relation
between the rational Casimir and KZ connections proved in \cite[Prop. 9.1]{TL2}.
\end{remark}

\subsection{The adjoint representation}

Drinfeld proved that, for any simple $\g$, the direct sum $\ol{\g}=\g\oplus\IC$
of the adjoint and trivial representations of $\g$ admits an extension to an
action of $Y(\g)$ on $\ol{\g}$ \cite[Thm. 8]{Dr2}. It is easy to check that
the corresponding action of $\H_\hbar'$ on $\ol{\g}[0]=\h\oplus\IC$ given
by Theorem \ref{th:trigo=AKZ} coincides with its action on affine linear
functions on $\h^*$ given by rational Dunkl operators (see, \eg \cite{Ki}).

\subsection{The case of $\sl{n}$}

Let $\g=\sl{n}$ and $V=\IC^n$ its vector representation. A simple inspection
shows that $V^{\otimes n}$ is a small \cite{Re2}. The zero weight space $V^{\otimes n}
[0]$ possesses two actions of the symmetric group: one arising from the Weyl
group action of $\SS_n$, the other from the permutation of the tensor factors,
under which it identifies with the group algebra $\IC\SS_n$.

The $\sl{n}$--module $V^{\otimes n}$ may be endowed with an action of $Y(\g)$
depending on $a_1,\ldots,a_n\in\IC$ obtained by composing the coproduct $\Delta
^{(n)}:Y(\g)\to Y(\g)^{\otimes n}$ with the evaluation homomorphisms $\ev_{a_i}:
Y(\g)\to U\g$. It is easy to check that the action of $\H'_\hbar$ on $V^{\otimes n}
[0]$ given by Theorem \ref{th:trigo=AKZ} coincides with that on the induced
representation $\ind^{\H'}_{S\h}\IC_{a_1,\ldots,a_n}$.

\section{Appendix: Tits extensions of affine Weyl groups}\label{se:Tits}

In this appendix, we review the definition of the Tits extension $\wt{W}$
of a Weyl group $W$. We then define the reduced Tits extension $\red
{W}$ of $W$ and show that, when $W$ is an affine Weyl group, $\red
{W}$ is isomorphic to the semi--direct product of the Tits extension of the
finite Weyl group underlying $W$ by the corresponding coroot lattice
(Theorem \ref{th:reduced affine Tits}).

\subsection{Weyl groups and braid groups}

Let $A=(a_{ij})_{i,j\in\bfI}$ be a generalised Cartan matrix and $(\h,
\Delta,\Delta^\vee)$ its unique realisation. Thus, $\h$ is a complex
vector space of dimension $2|\bfI|-\rk(A)$, $\Delta=\{\alpha_i\}_{i\in\bfI}
\subset\h^*$ and $\Delta^\vee=\{\cor{i}\}_{i\in\bfI}\subset\h$ are linearly
independent sets and
$$\<\cor{i},\alpha_j\>=a_{ij}$$
Recall that the Weyl group $W=W(A)$ attached to $A$ is the subgroup
of $GL(\h^*)$ generated by the reflections \cite[\S 3.7]{Ka}
$$s_i(\lambda)=\lambda-\<\lambda,\cor{i}\>\alpha_i$$
or, equivalently, the subgroup of $GL(\h)$ generated by the dual
reflections
$$s_i^\vee(t)=t-\<t,\alpha_i\>\cor{i}$$
By \cite[Prop. 3.13]{Ka}, the defining relations of $W$ are
\begin{align*}
s_i^2&=1\\
(s_i s_j)^{m_{ij}}&=1
\end{align*}
where for any $i\neq j$, $m_{ij}$ is equal to $2,3,4,6$ or $\infty$
according to whether $a_{ij}a_{ji}$ is equal to $0,1,2,3$ or $\geq
4$.

The {\it braid group} $B=B(A)$ attached to $A$ is the group with
generators $S_i$, $i\in\bfI$ and relations
$$\underbrace{S_i S_j\cdots}_{m_{ij}}=
\underbrace{S_j S_i\cdots}_{m_{ij}}$$
for any $i\neq j$.

\subsection{Tits extensions of Weyl groups}

\begin{definition}[\cite{Ti}]\label{de:Tits extension}
The {\it Tits extension} of $W$ is the group $\wt{W}$ with generators
$\wts_i$, $i\in\bfI$ and relations
\begin{align}
\underbrace{\wts_i\wts_j\cdots}_{m_{ij}}&=
\underbrace{\wts_j\wts_i\cdots}_{m_{ij}}
\label{eq:braid 1}\\
\wts_i^4&=1
\label{eq:braid 2}\\
\wts_i^2\wts_j^2&=\wts_j^2\wts_i^2
\label{eq:braid 3}\\
\wts_i\wts_j^2\wts_i^{-1}&=\wts_j^2(\wts_i^2)^{-a_{ji}}
\label{eq:braid 4}
\end{align}
\end{definition}

\subsection{}

Let $\g=\g(A)$ be the Kac--Moody algebra corresponding to the
Cartan matrix $A$ with generators $t\in\h$ and $e_i,f_i$, $i\in\bfI$.
Recall that a representation of $\g$ is integrable if $\h\subset\g$
acts semi--simply with finite--dimensional eigenspaces and 
$e_i,f_i$ act locally nilpotently. The next two results explain the
relevance and structure of the Tits extension $\wt{W}$.

\begin{proposition}\label{pr:Tits action}
Let $V$ be an integrable representation of $\g$. Then, the triple
exponentials
$$r_i=\exp(e_i)\exp(-f_i)\exp(e_i)$$
are well--defined elements of $GL(V)$ and the assignment
$\wts_i\rightarrow r_i$ yields a representation of $\wt{W}$
on $V$ mapping $\wts_i^2$ to $\exp(\pi\sqrt{-1}\cor{i})$.
\end{proposition}
\proof The $r_i$ are clearly well defined and satisfy
$$r_i \cdot t\cdot r_i^{-1}=s_i^\vee(t)$$
for any $t\in\h$ and $r_i^2=\exp(\pi\sqrt{-1}\cor{i})$
\cite[\S 3.8]{Ka}, from which \eqref{eq:braid 2}--\eqref{eq:braid 4} 
readily follow. Let now $i\neq j$ be such that $m_{ij}<\infty$. Then,
the Lie subalgebra $\g_{ij}$ of $\g$ generated by $e_i,f_i,\cor{i}$
and $e_j,f_j,\cor{j}$ is finite--dimensional and semi--simple and
$V$ integrates to a representation of the complex, connected and
simply--connected Lie group $G_{ij}$ with Lie algebra $\g_{ij}$.
By \cite{Ti}, $r_i$ and $r_j$ satisfy the braid relations \eqref{eq:braid 1}
when regarded as elements of $G_{ij}$\footnote{Tits' argument
is reproduced in the proof of (i) of Proposition \ref{pr:LG Tits}.},
and these therefore hold in $GL(V)$. \halmos

\subsection{}

Let $\Qv\subset\h$ be the lattice spanned by the coroots $\cor{i}$,
$i\in\bfI$. 

\begin{proposition}[\cite{Ti}]\label{pr:Tits extension}
$\wt{W}$ is an extension of $W$ by the abelian group $Z$
generated by the elements $\wts_i^2$. $Z$ is isomorphic,
as $W$--module to $\Qv/2\Qv\cong\IZ_2^{|\bfI|}$.
\end{proposition}
\proof Let $K\supset Z$ the kernel of the canonical projection
$\wt{W}\rightarrow W$. By \eqref{eq:braid 4}, $Z$ is a normal
subgroup of $\wt{W}$ and $\wt{W}/Z$ is generated by the images $\ol{s}_i$
of $\wts_i$ which, in addition to the braid relations \eqref{eq:braid 1},
satisfy $\ol{s}_i^2=1$. Thus, $\wt{W}/Z$ is a quotient of $W$,
$K=Z$ and $\wt{W}/Z\cong W$. Note next that, by \eqref
{eq:braid 2}--\eqref{eq:braid 4}, the assignment $\cor{i}
\rightarrow\wts_i^2$ extends to a $W$--equivariant surjection
$\Qv/2\Qv\rightarrow Z$. To prove that this is an
isomorphism it suffices to exhibit, for any $i\in\bfI$ a $\IZ_2$--valued
character $\chi_i$ of $Z$ such that $\chi_i(\wts_j^2)=(-1)^{\delta_{ij}}$.
Let $\lambda_i$ be the $i$th fundamental weight of $\g$, so that
$\<\lambda_i,\cor{j}\>=\delta_{ij}$, $V_i$ the irreducible
$\g$--module with highest weight $\lambda_i$ and $v_i\in V_i$
a nonzero highest weight vector. $V_i$ is integrable and since
$r_j^2=\exp(\sqrt{-1}\pi\cor{j})$, we have $r_j^2 v_i=(-1)^{\delta
_{ij}}v_i$. \halmos

\subsection{Reduced Tits extensions}

For any $v=\sum_i m_i\cor{i}\in \Qv$, set 
$$\wts^2_v=\prod_i(\wts_i^2)^{m_i}\in Z$$
so that for any $w\in W$ and lift $\wt{w}\in
\wt{W}$, $\wt{w}\wts_v^2\wt{w}^{-1}=\wts^2_{wv}$. By \cite
[Prop. 1.6]{Ka}, the center $\c$ of $\g$ is equal to
\begin{equation}\label{eq:center}
\c=\{t\in\h|\<t,\alpha_i\>=0\thickspace\text{for any $i\in\bfI$}\}
\end{equation}
The Weyl group operates trivially on $\c\subset\h$ and it follows
from \eqref{eq:braid 4} that the subgroup $Z_\c\subset Z$
generated by the elements $\wts^2_v$, with $v\in\c\cap \Qv
$ lies in the centre of $\wt{W}$.

\begin{definition}
The {\it reduced Tits extension} $\red{W}$ of $W$ is the quotient
$$\red{W}=\wt{W}/Z_\c$$
\end{definition}

By Proposition \ref{pr:Tits extension}, $\red{W}$ is an extension
of $W$ by $\Qv/(2\Qv+\c\cap \Qv)\cong\IZ_2^{\rk(A)}$.

\subsection{Reduced Tits extensions of affine Weyl groups}

Assume henceforth that $A=(a_{ij})_{0\leq i,j\leq n}$ is an affine
Cartan matrix of untwisted type. Altering our notations, we denote
by $\g$ the underlying complex, semi--simple Lie algebra and by
$\h$, $\{\alpha_i\}_{i=1}^{n}$, $\{\cor{i}\}_{i=1}^{n}$, $W$ and
$\Qv$ its Cartan subalgebra, simple roots, simple coroots, Weyl
group and coroot lattice respectively. Thus, for any $1\leq i,j\leq n$,
$$a_{ij}=\<\cor{i},\alpha_j\>,
\quad
a_{0j}=-\<\theta^\vee,\alpha_j\>
\quad\text{and}\quad
a_{j0}=-\<\cor{j},\theta\>$$
where $\theta\in\h^*$ is the highest root of $\g$. It is well known
that the (affine) Weyl group $W_a$ attached to $A$ is isomorphic
to the semi--direct product $W\ltimes \Qv$ \cite[prop. 6.5]{Ka}. The
isomorphism is given by mapping $s_i$ to $(s_i,0)$ for $i\geq 1$
and $s_0$ to $(s_\theta,-\theta^\vee)$.

The subspace $\c$ defined by \eqref{eq:center} is spanned by the
element
$$K=\cor{0}+\sum_{i=1}^n m_i\cor{i}$$
where the $m_i$ are the positive integers such that $\theta^\vee
=\sum_{i=1}^n m_i\cor{i}$ \cite[Prop. 6.2]{Ka}. It follows that the
reduced Tits extension $\red{\Wa}$ of $\Wa$ is the quotient of
$\wt{\Wa}$ by the relation
\begin{equation}\label{eq:braid 5}
\wts_0^2\cdot\prod_{i=1}^n (\wts_i^2)^{m_i}=1
\end{equation}

\subsection{Loop groups}\label{ss:start Wred}

The structure of the reduced Tits extension of $W_a$ will
be determined in \S \ref{ss:start Wred}--\S \ref{ss:end Wred}
by embedding $\red{\Wa}$ into the loop group corresponding
to $\g$.

Let $L\g=\g[z,z^{-1}]$ be the loop algebra of $\g$ and $d$ the
derivation of $L\g$ defined by $d x(m)=m x(m)$, where $x(m)
=x\otimes z^m$. Then, $L\g\rtimes\IC d$ is the quotient of the
Kac--Moody algebra corresponding to $A$ by the central element
$K$ defined above. Let $G$ be the complex, connected and
simply connected Lie group with Lie algebra $\g$ and $LG
=G(\IC[z,z^{-1}])$ the group of polynomial loops into $G$. Let
$H\subset G$ be the maximal torus with Lie algebra $\h$. The
group $\IC^*$ acts on $LG$ by reparametrisation fixing $G\supset
H$ and $H\times\IC^*$ is a maximal abelian subgroup of the
semi--direct product $LG\rtimes \IC^*$. By \cite{PS}, Prop. 5.2,
the normaliser of $H\times\IC^*$ in $LG\rtimes\IC^*$ is equal
to $(N(H)\ltimes H^\vee)\rtimes\IC^*$ where $N(H)$ is the
normaliser of $H$ in $G$ and $H^\vee=\Hom_\IZ(\IC^*,H)
\subset LG$ is isomorphic to the coroot lattice $\Qv$ by
$$\lambda\in \Qv\longrightarrow
\left(z\rightarrow z^{\lambda}=\exp_H(-\ln(z)\lambda)\right)$$
The quotient $N(H\times\IC^*)/H\times\IC^*$ is therefore isomorphic
to the affine Weyl group $\Wa=W\ltimes H^\vee$.

\subsection{}

For each real root $\wt{\alpha}=(\alpha,n)$ of $LG$, the subalgebra
$\sl{2}^{\wt{\alpha}}$ of $L\g$ spanned by
$$e_{\wt{\alpha}}=e_{\alpha}(n),\quad
f_{\wt{\alpha}}=f_{\alpha}(-n)\quad\text{and}\quad
h_{\alpha}$$
is the Lie algebra of a closed subgroup of $LG$ isomorphic to $SL_2(\IC)$.
This is obvious if $n=0$ and follows in the general case from the fact that
$\sl{2}^{(\alpha,n)}$ is conjugate to $\sl{2}^{(\alpha,0)}$. Indeed, any
element $\gamma_\lambda$ of the coweight lattice $\Hom(\IC^*,H/Z)
\subset L(G/Z)$ induces by conjugation an automorphism of $LG$
such that
$$\Ad(\gamma_\lambda)e_{\alpha}(n)=
e_{\alpha}(n-\<\lambda,\alpha\>)
\qquad\text{and}\qquad
\Ad(\gamma_\lambda)f_{\alpha}(n)=
f_{\alpha}(n+\<\lambda,\alpha\>)$$

\subsection{}

Let now $\alpha_i=(\alpha_i,0)$, $i=1\ldots n$ and $\alpha_0=(-\theta,1)$
be the simple roots of $LG$. For each $i=0\ldots n$, let $SL_2(\IC)\cong
G_i\subset LG$ be the corresponding subgroup, $H_i\subset G_i$ its
torus and $N_i$ the normaliser of $H_i$ in $G_i$. Note that any element
of $N_i\setminus H_i$ is of the form
$$\exp(e_i)\exp(-f_i)\exp(e_i)=\exp(-f_i)\exp(e_i)\exp(-f_i)$$
for some choice of root vectors $e_i\in(L\g)_{\alpha_i},f_i\in(L\g)_{-\alpha_i}$
such that $[e_i,f_i]=\cor{i}$ if $i\geq 1$ and $-\theta^\vee$ if $i=0$.

Let $\Ba$ be the (affine) braid group corresponding to $A$ and $S_0,S_1,
\ldots,S_n$ its generators.

\begin{proposition}\label{pr:LG Tits}\hfill
\begin{enumerate}
\item For any choice of $\sigma_i\in N_i\setminus
H_i$, $i=0\ldots n$, the assignment $S_i\rightarrow
\sigma_i$ extends uniquely to a homomorphism
$\sigma:\Ba\rightarrow N(H)\ltimes H^\vee$.
\item $\sigma$ factors through an isomorphism of the reduced
Tits extension $\red{\Wa}$ onto its image in $N(H)\ltimes H^\vee$.
\item If $\sigma,\sigma':\red{\Wa}\rightarrow N(H)\ltimes H^\vee$
are the homomorphisms corresponding to the choices $\{\sigma_i\}$
and $\{\sigma_i'\}$ respectively, there exists $t\in H\times\IC^*$
such that, for any $\wts\in\red{\Wa}$, $\sigma(\wts)=t\sigma'(\wts)
t^{-1}$.
\end{enumerate}
\end{proposition}
\proof
(1) the following argument is due to Tits \cite{Ti}. Let $i\neq j$ be
such that $m_{ij}$ is finite and set $s_{ij}=s_is_j\cdots\in\Wa$ and
$\sigma_{ij}=\sigma_i\sigma_j\cdots\in N(H)\ltimes H^\vee$ where
each product has $m_{ij}-1$ factors. The braid relations in $\Wa$
may be written as $s_{ij}s_{j'}=s_js_{ij}$ where $j'=j$ or $i$ according
to whether $m_{ij}$ is even or odd. Thus, $s_{ij}^{-1}s_js_{ij}=s_{j'}$
and therefore,
$$\Delta_{ij}=
\sigma_{j'}^{-1}\sigma_{ij}^{-1}\sigma_j\sigma_{ij}
\in H\cap\left(\sigma_{j'}^{-1}\sigma_{ij}^{-1}N_j\sigma_{ij}\right)
 = H\cap\sigma_{j'}^{-1}N_{j'}=H_{j'}$$
Repeating the argument with $i$ and $j$ permuted,
we find that $\Delta_{ji}\in H_{i'}$ with $i'=i$
or $j$ according to whether $m_{ij}$ is even or
odd. Thus, $\Delta_{ij}=\Delta_{ji}^{-1}\in H_{i'}
\cap H_{j'}=\{1\}$ where the latter assertion follows by follows from
the simple connectedness of $G$.

(2) The $\sigma_i$ satisfy \eqref{eq:braid 2}--\eqref{eq:braid 4}
and \eqref{eq:braid 5} since, for any $x_j\in N_j\setminus H_j$,
$x_j^2=\exp(i\pi\cor{j})$ for $j\geq 1$ and $x_0^2=\exp(-i\pi\theta
^\vee)$. Thus $\sigma$ descends to $\red{\Wa}$. Since the diagram
{\small $$\begin{xy}
\xymatrix@=0pt@C=0pt@R=12pt{
\red{\Wa}\ar[drrrrrr]\ar[rrrrrr]	&&&&&&	N(H)\ltimes H^\vee\ar[d]\\
						&&&&&&	\Wa}
\end{xy}$$}

\noindent
is commutative, the kernel of $\sigma$ is contained in $Z/Z_\c\cong
\IZ_2^n$ and is therefore trivial since, due to the simple--connectedness
of $G$, the subgroup of $G$ generated by $\sigma_j^2=\exp(\pi i\cor{j})$,
$j=1\ldots n$ is isomorphic to $\IZ_2^n$.

(3) For $i=1,\ldots,n$, let $t_i\in H_i$ be such that $\sigma_i=t_i\sigma_i'$
and choose $c_i\in\IC$ such that $t_i=\exp(c_ih_{\alpha_i})$. Since $s_i
\lambda_j^\vee=\lambda_j^\vee-\delta_{ij}\cor{i}$, where the $\lambda_j
^\vee\in\h$ are the fundamental coweights of $\g$, we find, with
$\ol{t}=\exp(\sum_{j=1}^n c_j\lambda_j^{\vee})\in H$,
$$\ol{t}\cdot\sigma_i'\cdot\ol{t}^{-1}
\exp(c_ih_{\alpha_i})\cdot\sigma_i'=
\sigma_i$$
Let now $t_0=\exp(c_0h_\theta)\in H_0$ be such that $\sigma_0=t_0\sigma
_0'$. Since for any $x\in\IC$,
$$\exp(xd)\sigma_0'\exp(-xd)=\exp(-xh_\theta)\sigma_0'$$
we find, with $y=\sum_j c_j\<\lambda_j^\vee,\theta\>-c_0$, that
$$\ol{t}\exp(yd)\sigma_0'\exp(-yd)\ol{t}^{-1}=
\exp((-y+\sum_jc_j\<\lambda_j^\vee,\theta\>)h_\theta)\sigma_0'=
\sigma_0$$
so that $t=\ol{t}\exp(yd)$ is the required element. \halmos

\subsection{}\label{ss:end Wred}

The following is the main result of this appendix.

\begin{theorem}\label{th:reduced affine Tits}
The inclusion $\wt{W}\hookrightarrow\red{W_a}$ extends to an isomorphism
$\wt{W}\ltimes\Qv\to\red{W_a}$ making the following a commutative diagram
{\small $$\begin{xy}
\xymatrix@=0pt@C=20pt@R=15pt{
\wt{W}\ltimes\Qv\ar[r]\ar[d] 	& \ar[d]\red{W_a}\\
W\ltimes\Qv\ar@{=}[r]		& W_a}
\end{xy}$$}
\end{theorem}
\proof We wish to construct a $\wt{W}$--equivariant section $s$ to the
restriction to $\Qv$ of the extension 
$$1\rightarrow Z/Z_c\rightarrow\red{\Wa}\rightarrow\Wa\rightarrow 1$$
Identify for this purpose $\red{\Wa}$ with its image inside $N(H\times\IC
^*)\cap LG$ by using Proposition \ref{pr:LG Tits}. We claim that there
exists $x\in\IC$ such that, for any $\lambda\in \Qv$, $\exp_H(x\lambda)
\cdot z^\lambda$ lies in $\red{\Wa}$. It is then clear that $s(\lambda)=
\exp_H(x\lambda)\cdot z^\lambda$ yields the required section.

Let $\cor{i}$ be a short simple coroot and $w\in W$ an element such
that $\theta^\vee=w\cor{i}$. Let $\wt{w}\in\wt{W}$ be a lift of
$w$ and lift $\theta^\vee\in \Qv$ to
$$\tau^{\theta^\vee}=\wts_0\,\wt{w}\wts_i\wt{w}^{-1}\in\red{\Wa}$$
Let $e_\theta\in\g_\theta,f_\theta\in\g_{-\theta}$ be root vectors such
that $[e_\theta,f_\theta]=\theta^\vee$ and denote by $\rho_{\theta^{\vee}}:
SL_2\to G$ the embedding whose differential maps $e,f,h\in\sl{2}$
to $e_\theta,f_\theta,\theta^\vee$. We may assume $e_\theta,f_\theta$
chosen so that $\wts_0$ is of the form
$$\exp(f_\theta\otimes z)\exp(-e_\theta\otimes z^{-1})\exp(f_\theta\otimes z)
=\rho_{\theta^\vee}(
\begin{pmatrix}0& -z^{-1}\\z&0 \end{pmatrix})$$
Since
$\wt{w}\wts_i\wt{w}^{-1}\in N_\theta\setminus H_\theta$
is necessarily of the form 
$$\exp(te_\theta)\exp(-t^{-1}f_\theta)\exp(te_\theta)=
\rho_{\theta^\vee}(\begin{pmatrix}0& t\\-t^{-1}&0 \end{pmatrix})$$
for some $t\in\IC^*$, we find that
$$\tau^{\theta^\vee}=
\rho_{\theta^\vee}(\begin{pmatrix}(tz)^{-1}&0\\0&tz\end{pmatrix})=
\exp(x\theta^\vee)\cdot z^{\theta^\vee}$$
with $x=-\ln(t)$ which proves our claim for $\lambda=\theta^\vee$.
Let now $w\in W$ with lift $\wt{w}\in\wt{W}$, then
$$\wt{w}\tau_{\theta^\vee}\wt{w}^{-1}=\exp(xw(\theta^\vee))\cdot
z^{w\theta^\vee}$$
so that that $\exp(x\cor{i})\cdot z^{\cor{i}}\in\red{\Wa}$ for any
short coroot $\cor{i}$. Since the short coroots span $\Qv$, the
claim holds for any $\lambda\in \Qv$.

Since $s(\Qv)$ is free abelian and $\wt{W}$ is finite, their intersection
is trivial and the map $\wt{W}\ltimes \Qv\rightarrow\red{\Wa}$, $(\wt{w},
\lambda)\rightarrow \wt{w}s(\lambda)$ is injective. It is moreover
surjective since $Z/Z_c$ is generated by $\wts_i^2$, $i=1\ldots n$
and therefore lies in $\wt{W}$ \halmos

\subsection{}

\begin{remark}
Unlike $\red{\Wa}$\negthickspace\negthinspace, the (non--reduced)
Tits extension $\wt{\Wa}$ of $\Wa$ is not a semi--direct product
in general. For example, for $\g=\sl{2}$, with affine Cartan matrix
$$A=\left(\begin{array}{rr}2&-2\\-2&2\end{array}\right)$$
$\wt{\Wa}$ is generated by $\wts_0$, $\wts_1$ with relations
$\wts_i^4=1$,
$$\wts_0\wts_1^2\wts_0^{-1}=\wts_1^2(\wts_0^2)^2=\wts_1^2
\quad\text{and}\quad
\wts_1\wts_0^2\wts_1^{-1}=\wts_0^2$$
In particular, the group $Z\cong\IZ_2^2$ generated by $\wts_0^2,
\wts_1^2$ lies in the centre of $\wt{\Wa}$. Any lift in $\wt{\Wa}$ of
the generator of $\Qv\cong\IZ$ is of the form $\tau=z\wts_0\wts_1$,
for some $z\in Z$ and gives rise a $\wt{W}$--equivariant section
$\Qv\rightarrow\wt{\Wa}$ if, and only if, $\wts_1\tau\wts_1^{-1}=
\tau^{-1}$. Since $z=z^{-1}$ is central, such a section exists iff
$\wts_1(\wts_0\wts_1)\wts_1^{-1}=\wts_1^{-1}\wts_{0}^{-1}$ and
therefore iff $\wts_1^2\wts_0^{2}=1$ which holds in $\red{\Wa}$
but not in $\wt{\Wa}$.
\end{remark}

\begin{remark} The section $\Qv\rightarrow\red{\Wa}$ constructed
in Theorem \ref{th:reduced affine Tits} does not in general coincide
with that obtained from the canonical section $\Qv\rightarrow\Ba$
\cite[\S 3.2--3.3]{Mc}. For example, for $\g=\sl{3}$, the canonical lift of $\theta
^\vee\in \Qv$ in $\Ba$ is $T^{\theta^\vee}=S_0 S_1S_2S_1$. When
regarded as an element $\tau^{\theta^\vee}$ of $\red{\Wa}$ this
does not give rise to a $\wt{W}$--equivariant section since
$\Ad(\wts_\theta)\tau^{\theta^\vee}\neq (\tau^{\theta^\vee})^{-1}$
where $\wt{s}_\theta=\wts_1\wts_2\wts_1$ is a lift in $\wt{W}$ of
the reflection $s_\theta$. Indeed, $\Ad(\wts_2)\wts_1^2=\wts _1
^2\wts_2^2$ in $\wt{W}$, so that
$$\wt{s}_\theta^2=
\wts_1\Ad(\wts_2)(\wts_1^2)\wts_2^2\wts_1=
\wts_1^3\wts_2^4\wts_1=
1$$
Thus, since $\tau^{\theta^\vee}=\wts_0\wts_\theta$,
$$\Ad(\wts_\theta)\tau^{\theta^\vee}=\wts_\theta\wts_0
\qquad\text{while}\qquad
(\tau^{\theta^\vee})^{-1}=\wts_\theta\wts_0^{-1}$$
which are different elements of $\red{\Wa}$ by \eqref
{eq:braid 5}. 
\end{remark}

\section*{Acknowledgments}

I am grateful to Rapha\"{e}l Rouquier for pointing me towards the
Yangian at the early stages of this project. I am also grateful to
Sachin Gautam for spotting a gap in my initial proof of Theorem
\ref{th:trigo flat} and to the referee for a very careful reading of this
manuscript. This paper was begun at the summer home of P\'en\'elope
Riboud and Romain Graziani. It is a pleasure to thank them for their
warm and friendly hospitality, and their daughter L\'eonore for a
number of inspiring comments.

\end{document}